\documentclass{amsart}
\newcommand{\note}{\noindent {\bf Notation. }}

\newcommand{\ws}{\hspace{4pt}}
\newtheorem{theorem}{Theorem}

\newtheorem{remark}{Remark}
\newtheorem{proposition}{Proposition}
\newtheorem{cor}{Corollary}
\newtheorem{lemma}{Lemma}

\newtheorem{defi}{Definition}
\usepackage[]{amssymb,amsopn}
\usepackage[]{amsmath}
\usepackage[]{eucal}
\usepackage[]{latexsym}
\usepackage{color}

\makeatletter
\@namedef{subjclassname@2020}{%
  \textup{2020} Mathematics Subject Classification}
\makeatother

\begin{document}

\title[B-$p$ capacity]{$p$-capacity with Bessel convolution }
\author{\'A. P. Horv\'ath}

\subjclass[2020]{31C45, 26D15, 28A78}
\keywords{nonlinear potential, Bessel convolution, Laplace-Bessel equation, Wolff inequality, weighted Hausdorff measure}
\thanks{}

\begin{abstract}
We define and examine nonlinear potential by Bessel convolution with Bessel kernel. We investigate removable sets with respect to Laplace-Bessel inequality. By studying the maximal and fractional maximal measure, a Wolff type inequality is proved. Finally the relation of B-$p$ capacity and B-Lipschitz mapping, and the B-$p$ capacity and weighted Hausdorff measure and the B-$p$ capacity of Cantor sets are examined.
\end{abstract}
\maketitle

\section{Introduction}

Classical, nonlinear, and Bessel potentials are widespread, have an extensive literature, and are widely applicable, see e.g. \cite{rs}, \cite{ghn}, \cite{mp} and the references therein. Below we introduce and examine nonlinear potential defined by Bessel convolution with Bessel kernel.

Bessel translation was defined by Delsarte \cite{d} and the basic investigation is due to Levitan, \cite{le}. In a series of works the authors pointed out that Bessel translation and convolution methods are effective tools  to handle Bessel-type partial differential operators, see e.g. \cite{p}, \cite{kk}, \cite{esk}, \cite{ly}. It also proved useful for deriving Nikol'skii type inequality, see \cite{abdh}, and for giving compactness criteria in some Banach spaces, see \cite{h}.

This leads to examine nonlinear potential and $p$-capacity with respect to Bessel convolution. The curiosity of the method is that the underlying space of Bessel-$p$ capacity is automatically weighted. Weighted nonlinear potential was studied already in the '80-s, see e. g. \cite{a}, \cite{aclm}. For logarithmic potentials with external field see the monograph \cite{sat}. In our investigation the Bessel weighted space is a natural consequence of the definition of convolution, and so many of the results are very similar to the ones proved in the unweighted case.

The paper is organized as follows. After the preliminaries, in the third section, applying recent results on Bessel potential, we investigate removable sets for Laplace-Bessel equation. In the fourth section a Wolff type inequality is proved, which is the basis of the study of the last section. This last section contains some "metric" results on Lipschitz type mapping and on capacity of Cantor sets. Since Bessel translation is not a geometric similarity, moreover the underlying space is weighted, we have to introduce a special property (B-Lipschitz mapping), and the notion of weighted Hausdorff measure.

\section{Notation, preliminaries}
Let $\mathbb{R}^n_+:=\{ x=(x_1,\dots ,x_n): x_i\ge 0, \ws i=1, \dots , n.\}$. $\lambda$ is the $n$-dimensional Lebesgue measure. $a=a_1,\dots ,a_n$ is a multiindex. Let $E\subset \mathbb{R}^n_+$ and $\mathcal{M}(E)$ stands for the Radon measures supported on $E$. If $\mu \in \mathcal{M}(E)$ for some $E$, $d\mu_a(x):=x^ad\mu(x)$
Define the Banach space $L^p_a$ as follows.
$$\|f\|_{p,a}^p=\int_{\mathbb{R}^n_+}|f(x)|^pd\lambda_a(x),$$
and as usual
$$L^p_a:=L^p_a(\mathbb{R}^n_+)=\{f : \|f\|_{p,a}<\infty\}, \ws \ws L^{p+}_a:=\{f\in L^p_a : f\ge 0\}.$$
The dual index $p'$ is defined by $\frac{1}{p}+\frac{1}{p'}=1$.

\medskip

\subsection{Bessel translation}
Let
$$a:=2\alpha_1+1, \dots , 2\alpha_n+1,  \ws \ws \ws \alpha_i>-\frac{1}{2}, \ws i=1, \dots , n, \ws \ws |a|=\sum_{i=1}^n(2\alpha_i+1).$$
The Bessel translation of a function, $f$ (see e.g. \cite{le}, \cite{p}, \cite{ss}) is
$$T_{a}^tf(x)=T_{a_n}^{t_n}\dots T_{a_1}^{t_1}f(x_1,\dots ,x_n),$$
where
\begin{equation}\label{tra1}T_{a_i}^{t_i}f(x_1,\dots ,x_n)$$ $$=\frac{\Gamma(\alpha_i+1)}{\sqrt{\pi}\Gamma\left(\alpha_i+\frac{1}{2}\right)}\int_0^\pi f(x_1,\dots,\sqrt{x_i^2+t_i^2 -2x_it_i\cos\vartheta_i}, x_{i+1}, \dots, x_n)\sin^{2\alpha_i}\vartheta_i d\vartheta_i.\end{equation}
The translation can also be expressed as an integral with respect to a kernel function:
\begin{equation}\label{tra2}T_{a_i}^{t_i}f(x_1,\dots ,x_n)=\int_0^\infty K(x_i,t_i,z_i)f(z_1,\dots ,z_n)d\lambda_{a_i}(z_i),\end{equation}
where
\begin{equation}\label{trake}K(x,t,z)=\left\{\begin{array}{ll}\frac{\pi^{\alpha+\frac{1}{2}}\Gamma(\alpha+1)}{2^{2\alpha-1}\Gamma\left(\alpha+\frac{1}{2}\right)}\frac{[(x+t)^2-z^2)(z^2-(x-t)^2)]^{\alpha-\frac{1}{2}}}{(xtz)^{2\alpha}}, \ws\ws |x-t|<z<x+t\\
0, \ws \ws \mbox{otherwise}. \end{array}\right.\end{equation}

Obviously
$$T_{a}^tf(x)=T_{a}^xf(t).$$
$T_a$ is a positive operator, and
\begin{equation}\label{t1}\|T_{a,x}^tf(x)\|_{p,a}\le \|f\|_{p,a}, \ws \ws 1\le p \le\infty,\end{equation}
see e.g. \cite{le}.

\medskip

The generalized convolution with respect to the Bessel translation is
$$f*_a g =\int_{\mathbb{R}^n_+}T_{a,x}^tf(x)g(x)d\lambda_a(x).$$
We have
$$f*_a g =g*_a f,$$
and Young's inequality fulfils i. e. if $1\le p, q, r \le \infty$ with $\frac{1}{r}=\frac{1}{p}+\frac{1}{q}-1$; if $f\in L^p_\alpha$ and $g \in L^q_\alpha$, then
\begin{equation}\label{Y}\|f*_ag\|_{r,a} \le \|f\|_{p,a}\|g\|_{q,a},\end{equation}
see \cite[(3.178)]{ss}.

\medskip

Subsequently if it does not cause any confusion, $T^tf(x)$ stands for $T_{a}^tf(x)$. For any set $H\subset\mathbb{R}^n$ we denote by $H_+:=H\cap \mathbb{R}^n_+$. The next technical lemma will be useful in the following sections.

\medskip

\begin{lemma}\label{tk} $\mathrm{supp}T^t\chi_{B_+(0,r)}(x)=\overline{B_+(x,r)}$, $\mathrm{supp}T^t\chi_{[0,r)^n}(x)=\times_{i=1}^n[x_i-r,x_i+r]_+=:T_+(x,r)$.
There is a $c>0$ such that for all $x\in \mathbb{R}^n_+$, $t\in B_+(x,r)$
\begin{equation}\label{kf}T^t\chi_{B_+(0,r)}(x)\le c \prod_{i=1}^n\min\left\{1,\left(\frac{r}{x_i}\right)^{a_i}\right\}.\end{equation}
There is a $c>0$ such that for all $x\in \mathbb{R}^n_+$, $t\in T_+\left(x,\frac{r}{2}\right)$
\begin{equation}\label{kf1}T^t\chi_{[0,r)^n}(x)\ge c \prod_{i=1}^n\min\left\{1,\left(\frac{r}{x_i}\right)^{a_i}\right\}.\end{equation}
\end{lemma}

\proof The first two statements are direct consequences of the definition, for \eqref{kf} see \cite[p. 321]{g}.\\
\eqref{kf1}: Since $\sqrt{x_i^2+t_i^2-2x_it_i\cos\vartheta_i}\le x_i+t_i$ if $x_i+t_i\le r$, then
$$\int_{\{\vartheta \in [0\pi) : \sqrt{x_i^2+t_i^2-2x_it_i\cos\vartheta_i}\le r\}}1d\sigma_i\vartheta=1,$$
where $d\sigma_i:=\frac{\Gamma(\alpha_i+1)}{\sqrt{\pi}\Gamma\left(\alpha_i+\frac{1}{2}\right)}\sin\vartheta^{2\alpha_i}d\vartheta_i$ is a probability measure on $[0,\pi]$.\\
If $r<x_i+t_i\le 2r$, using \eqref{tra2} and recalling that $|x_i-t_i|\le \frac{r}{2}$ we have
$$I_i=\frac{c(a_i)}{(x_it_i)^{2\alpha_i}}\int_{|x_i-t_i|}^r\left[(z_i^2-(x_i-t_i)^2)((x_i+t_i)^2-z_i^2)\right]^{\alpha_i-\frac{1}{2}}z_idz_i$$ $$\ge c(a_i)\frac{1}{(x_i+t_i)^{4\alpha_i}}\int_{\frac{5}{8}r}^{\frac{7}{8}r}(\cdot)dz_i= c \frac{1}{r^{4\alpha_i}}r^{4\alpha_i-2}r r\ge c(a_i).$$
If $x_i+t_i> 2r$, then $x_i\sim t_i$ ($f\sim g$ if there are positive constants $A$ and $B$ such that $Af<g <Bf$) and we have
$$I_i\ge c \int_{\frac{5}{8}r}^rr^{2\alpha_i-1}\frac{(x_i+t_i)^{2\alpha_i-1}}{(x_it_i)^{2\alpha_i}}z_idz_i =c\left(\frac{r}{x_i}\right)^{2\alpha_i+1}.$$

\medskip

\begin{remark}\label{cr} {\rm For any $c\ge 2\sqrt{n}$
$$T^t\chi_{B_+(0,r)}(x)\sim T^t\chi_{B_+(0,cr)}(x),$$
for all $x\in \mathbb{R}^n_+$, $t\in B_+(x,r)$.\\
Indeed, in view of \eqref{tra2} if $H\subset S\subset \mathbb{R}^n_+$, then $T^t\chi_H(x)\le T^t\chi_S(x)$ for any $x,t$. Together with Lemma \ref{tk} it implies that $c_1T^t\chi_{B_+(0,r)}(x)\le T^t\chi_{B_+(0,cr)}(x)\le c_2 T^t\chi_{B_+(0,r)}(x)$.}\end{remark}

\medskip

\subsection{Radially decreasing kernels and B-$p$ capacity}

\begin{defi}\label{def1}Let $g$ be a non-negative lower semi-continuous, non-increasing function on $\mathbb{R}_+$ for which
\begin{equation}\label{i0}\int_0^1g(t)t^{n+|a|-1}dt<\infty.\end{equation}
Then $\kappa:=g(|x|)$ is a radially decreasing kernel on $\mathbb{R}^n$.
\end{defi}

\medskip

The B-$p$ capacity with respect to $\kappa$ is as follows.

\begin{defi}\label{cap} Let $E\subset \mathbb{R}^n_+$, $1\le p <\infty$.
$$C_{p,\kappa}(E):=\inf\{\|f\|_{p,a}^p: f\in L^{p +}_a, \ws \kappa *_a f(x)\ge 1,  \ws \forall x\in E \}.$$
\end{defi}

\medskip

\begin{remark}\label{cp}
{\rm (1) Definition \ref{def1} is a special case of \cite[Definition 2.3.3]{ah}. Thus all the standard properties proved in \cite[Chapter 2.3]{ah} are valid.\\
(2) Notice that by the definition of Bessel translation if $\kappa$ is a radially decreasing kernel, then $T^t\kappa(x)\le g(|x-t|)$. Thus if $f\ge 0$, $\kappa *_a f(x)\le\kappa*h(x)$, where $h(x)=f(x)x^a\chi_{\mathbb{R}^n_+}(x)$ and $*$ stands for the standard convolution.\\
(3) Let $K\subset \mathbb{R}^n_+$, $1<p<\infty$. An equivalent form of Definition \ref{cap} is
$$C_{p,\kappa}^{\frac{1}{p}}(K)=\sup\{\mu_a(K) : \mu \in \mathcal{M}(K), \ws \|\kappa*_a\mu\|_{p',a}\le 1\},$$
where $\mathcal{M}(K)$ is the set of (positive) measures on $K$, see \cite[Theorem 2.5.1]{ah}.\\
(4) As usual, the definitions above can be extended to any subsets of $\mathbb{R}^n_+$ as it follows. If $O\subset \mathbb{R}^n_+$ is open, then $C_{p,\kappa}(O):=\sup \{C_{p,\kappa}(K): K\subset O, \ws K \ws \mbox{is compact}\}$ and if $E\subset \mathbb{R}^n_+$ is arbitrary, then $C_{p,\kappa}(E):=\inf \{C_{p,\kappa}(O): E\subset O, \ws O \ws \mbox{is open}\}$.\\
(5) $C_{p,\kappa}$ is monotone and $\sigma$-subadditive (cf. \cite[Propositions 2.3.4 and 2.3.6]{ah}).\\}
\end{remark}

\medskip

\begin{proposition}\label{pro} Let $1<p<\infty$. If $\kappa$ is a radially decreasing kernel, then\\
\rm{(1)} if $\|\kappa\|_{p',a}<\infty$, then $C_{p,\kappa}(\{y\})>0$ for all $y\in \mathrm{int}\mathbb{R}^n_+$,\\
\rm{(2)} if $\int_{\mathbb{R}^n_+ \setminus B(0,1)}\kappa^{p'}d\lambda_a =\infty$, then $C_{p,\kappa}(E)=0$ for all $E\subset \mathbb{R}^n_+$,\\
\rm{(3)} if $\int_{\mathbb{R}^n_+ \setminus B(0,1)}\kappa^{p'}d\lambda_a <\infty$, $E$ is measurable and $C_{p,\kappa}(E)=0$, then $\lambda_a(E)=0$.
\end{proposition}

\proof \rm{(1)} Let $\delta_y$ be the Dirac measure concentrated at $y$.  According to Remark \ref{cp} and \eqref{Y},
$$C_{p,\kappa}^{\frac{1}{p}}(\{y\})=\sup\left\{\frac{\mu_a(\{y\})}{\|\kappa*_a \delta_y\|_{p',a}} :  \mu \in \mathcal{M}(\{y\})\right\}\ge \frac{y^a}{\|\kappa*_a \delta_y\|_{p',a}}\ge \frac{1}{\|\kappa\|_{p',a}}>0.$$
\rm{(2)} It is enough to show that $C_{p,\kappa}(B_+(0,r))=0$ for all $r>0$. Let $\mu\in\mathcal{M}(B_+(0,r))$. In view of \eqref{tra1} $T^t\kappa(x)\ge g(|x+t|)$, thus
$$\|\kappa *_a \mu\|_{p',a}\ge \mu(B_+(0,r))\left(\int_{\mathbb{R}^n_+}g(r+|t|)^{p'}t^adt\right)^{\frac{1}{p'}}$$ $$\ge c \left(\int_{\mathbb{R}^n_+ \setminus B(0,2r)}g(|t|)^{p'}t^adt\right)^{\frac{1}{p'}}=\infty.$$
The last inequality is equivalent with the assumption, and according to Remark \ref{cp}, it proves the statement.\\
\rm{(3)} It is enough to show that $\lambda_a(E\cap B_+(0,r))=0$ for all $r>0$. Let $F=E\cap B_+(0,r)$ and $f\in L^{p +}_a$ such that $\kappa *_a f(x)\ge 1$ on $F$. Then by Fubini's theorem
$$\lambda_a(F)\le\int_F \kappa *_a f(x)x^adx=\int_{\mathbb{R}^n_+}\chi_F(x)\kappa *_a f(x)x^adx=\int_{\mathbb{R}^n_+}\kappa *_a\chi_F(t)f(t)t^adt$$ $$\le \|f\|_{p,a}\|\kappa *_a\chi_F\|_{p',a}\le  \|f\|_{p,a}\|\kappa *_a\chi_{B_+(0,r)}\|_{p',a}.$$
We estimate the second factor.
$$\|\kappa *_a\chi_{B_+(0,r)}\|_{p',a}$$ $$\le  \left(\int_{B_+(0,2r)}(\kappa *_a\chi_{B_+(0,r)}(t))^{p'}t^adt\right)^{\frac{1}{p'}}+\left(\int_{\mathbb{R}^n_+ \setminus B_+(0,2r)}(\cdot)\right)^{\frac{1}{p'}}=I+II.$$
If $|t|>2r$, $\kappa(x)\le\kappa\left(\frac{t}{2}\right)$ while by Lemma \ref{tk} $T^t\chi_{B_+(0,r)}(x)\le c r^{|a|}\frac{1}{x^a}$ on $|x-t|<r$. Thus by the assumption we have
$$II=\left(\int_{\mathbb{R}^n_+ \setminus B_+(0,2r)}\left(\int_{\mathbb{R}^n_+}T^t\chi_{B_+(0,r)}(x)\kappa(x)x^adx\right)^{p'}t^adt\right)^{\frac{1}{p'}}$$ $$\le c\left(\int_{\mathbb{R}^n_+ \setminus B_+(0,2r)}\kappa^{p'}\left(\frac{t}{2}\right)t^a\right)^{\frac{1}{p'}}< c.$$
In the first integral $|t|<2r$ and $|x-t|<r$, so the convolution can be estimated as
$$\kappa *_a\chi_{B_+(0,r)}(t)\le \int_{B_(0,3r)}(g|x|)x^adx\le c,$$
where in spherical coordinates the last inequality is just \eqref{i0}. Thus, $I$ is also bounded by a constant. Taking infimum over appropriate functions $f$, we have that $\lambda_a(F)\le c C_{p,\kappa}(F)$, which implies the statement.

\medskip

\begin{remark}{\rm Of course, the nonlinear potential with Bessel convolution is $V_{\kappa,p}^\mu=\kappa*_a(\kappa*_a\mu)^{p'-1}$. Subsequently, we focus on capacity.
}\end{remark}

\medskip

\subsection{Bessel and Riesz kernels}

The modified Bessel function of the second kind, $K_\alpha$ is defined as follows.
$$i^{-\alpha}J_\alpha(ix)=\sum_{k=0}^\infty \frac{1}{k!\Gamma(k+\alpha+1)}\left(\frac{x}{2}\right)^{2k+\alpha},$$
where $J_\alpha$ is the Bessel function, and
$$K_\alpha(x)=\frac{\pi}{2}\frac{i^{\alpha}J_{-\alpha}(ix)-i^{-\alpha}J_\alpha(ix)}{\sin\alpha \pi}.$$
Considering $r>0$, around zero
\begin{equation}\label{k1}K_\alpha(r)\sim \left\{\begin{array}{ll}-\ln\frac{r}{2}-c, \ws \mbox{if}\ws  \alpha=0\\ C(\alpha)r^{-\alpha},  \ws \mbox{if}\ws \alpha>0,\end{array}\right.\end{equation}
and around infinity
\begin{equation}\label{k2}K_\alpha(r)\sim \frac{c}{\sqrt{r}}e^{-r}.\end{equation}
The Bessel kernel is
\begin{equation}\label{g1}G_{a,\nu}(x):=\frac{2^{\frac{n-a-\nu}{2}+1}}{\Gamma\left(\frac{\nu}{2}\right)\prod_{i=1}^n\Gamma(\alpha_i+1)}\frac{K_{\frac{n+|a|-\nu}{2}}(|x|)}{|x|^{\frac{n+|a|-\nu}{2}}}.\end{equation}

\medskip

Below we also need the Riesz kernel:
\begin{equation}I_{\beta}(x)=\frac{c(\beta)}{|x|^{n-\beta}}, \ws \ws x\in \mathbb{R}^n.\end{equation}
In the last section we use Bessel kernel rather than the Riesz kernel, because its behavior at infinity allows wider function classes. On the other hand, around the origin the Riesz kernel, $I_{\nu-|a|}(x)$, behaves similarly to the Bessel kernel and is simpler, thus it proved to be a useful tool for computations.

\medskip

Below we examine B-$p$ capacity, which is defined by generalized convolution referring to the Bessel kernel: $C_{p,a,\nu}(E):=C_{p,G_{a,\nu}}(E)$.
In view of \eqref{k1}, \eqref{k2} and \eqref{g1} $G_{a,\nu}(x) \in L^1_a$ if and only if $\nu>0$. On the other hand according to Proposition \ref{pro}, B-$p$ capacity is non trivial if and only if $1<p< \frac{n+|a|}{\nu}$ or $1=p=\frac{n+|a|}{\nu}$. Thus subsequently we investigate $C_{p,a,\nu}$ if $1<p<\infty$, that is
\begin{equation}\label{g2}0<\nu< \frac{n+|a|}{p}.\end{equation}

\medskip

\section{The Laplace-Bessel operator}

B-elliptic equations are investigated by several authors. For instance fundamental solutions are given, see e.g.\cite{kk} and \cite{ly}. Harmonic analysis associated with Bessel operator is examined, see e.g. \cite{p} and mean-value theorems are proved, see \cite{s}. Here we give a simple application of B-$p$ capacity.

We begin this section by introducing some additional notation.
According to  \eqref{g2}, \eqref{t1}   if $g\in L^p_a$, $G_{a,\nu}*_a g\in L^p_a$, moreover by \cite[Lemma 4.3 (3)]{esk}
\begin{equation} \label{gn}\|G_{a,\nu}*_a g\|_{p,a}\le\|g\|_{p,a}.\end{equation}
Thus we define the next Banach space.
$$L^p_{a,\nu}:=L^p_{a,\nu}(\mathbb{R}^n_+)=\{f :f=G_{a,\nu}*_a g; \ws   g\in L^p_a\}, \ws \ws \ws \|f\|_{p,a,\nu}:=\|g\|_{p,a}.$$
Let
$$B_{\alpha,x}:=\frac{\partial^2}{\partial x^2}+\frac{2\alpha+1}{x}\frac{\partial}{\partial x}$$
be The Bessel operator. The Laplace-Bessel operator is defined as
$$\Delta_a=\sum_{i=1}^n B_{\alpha_i, x_i}.$$
With this notation we define the Sobolev space $W^m_{p,a}$ with $m \in\mathbb{N}$ as it follows.
$$W^m_{p,a}:=\{f\in L^p_a : \Delta_a^k f \in  L^p_a, \ws k=1,\dots, m\}, \ws \ws \ws \|f\|_{W^m_{p,a}}=\sum_{k=0}^m \|\Delta_a^k f\|_{p,a}.$$

\medskip

\note
We need the "even" functions from the Schwartz class in $\mathbb{R}^n_+$.
$$\mathcal{S}_e:=\left\{f\in C^\infty(\mathbb{R}^n_+): \left.\frac{\partial^{2k+1}f}{\partial x_i^{2k+1}}\right|_{x_i=0}=0, \ws k\in\mathbb{N};\right.$$ $$\left.    \sup_{x\in\mathbb{R}^n_+}\left|x^\alpha D^{\beta}f(x)\right|<\infty \ws \forall \alpha, \beta \in \mathbb{N}^n\right\}.$$

\medskip

The next lemma describes the relation of Bessel potential and Sobolev spaces above.

\begin{lemma}\label{l1} Let $m$ be a positive integer. Then
$$W^m_{p,a}=L^p_{a,2m}.$$
\end{lemma}

\proof Let us notice first that if we define $\|f\|_{\tilde{W}^m_{p,a}}=\sum_{k=0}^m \|(I-\Delta_a)^k f\|_{p,a}$, then $\|f\|_{\tilde{W}^m_{p,a}}\sim \|f\|_{W^m_{p,a}}$. According to \cite[Lemma 4.3]{esk} if $g \in  L^p_a$, $1\le p \le \infty$ and $k\in \mathbb{N}$,
\begin{equation}\label{s1} (I-\Delta_a)^k(G_{a,\nu}*_a g)=G_{a,\nu-2k}*_a g, \ws \ws \mbox{and} \ws G_{a,0}*_a g=g.\end{equation}
Comparing this with \eqref{gn} we have if $f\in L^p_{a,2m}$ and $k\le 2m$, then
$$\|(I-\Delta_a)^kf\|_{p,a}\le \|G_{a,2m-2k}*_a g\|_{p,a}\le \|g\|_{p,a}= \|f\|_{p,a,2m}.$$
On the other hand taking into consideration that  $S_e$ is a dense subset in $W^m_{p,a}$, let $f\in S_e$. According to \cite[Theorem 4.5]{esk} for $1\le p \le\infty$, $k\in \mathbb{N}$
\begin{equation}\label{s2}G_{a,\nu+2k}*_a (I-\Delta_a)^k f =G_{a,\nu}*_a f.\end{equation}
Thus by \eqref{s1} and \eqref{s2}
$$f= G_{a,0}*_a f=G_{a,2m}*_a(I-\Delta_a)^mf=G_{a,2m}*_a g,$$
where $g=(I-\Delta_a)^mf \in  L^p_a$. So $f\in L^p_{a,2m}$, and
$$\|f\|_{p,a,2m}\le \|f\|_{W^m_{p,a}}.$$

\medskip

\begin{defi}\label{N}Let $K\subset \mathbb{R}^n_+$ compact. $\mathcal{S}$ is the Schwartz class restricted to $\mathbb{R}^n_+$.
$$N_{p,a,\nu}(K):=\inf\{\|f\|_{p,a,\nu}^p: f=G_{a,\nu}*_a g\in\mathcal{S}, \ws f\equiv 1 \ws \mbox{in a neigborhood of}  \ws K\}.$$
\end{defi}

\begin{remark}
{\rm (1) If $p=2$ with standard convolution, $N$ is the spectral measure defined by Deny, see \cite{de}.\\
(2)Of course, $N_{p,a,\nu}$ can be extended as above and $C_{p,a,\nu}(E)\le N_{p,a,\nu}(E)$.}
\end{remark}

\medskip

\note
 Let us introduce the inner product for measurable functions
 $$\langle f,g\rangle_a := \int_{\mathbb{R}^n_+} fgd\lambda_a.$$
 We denote in the same way the effect of a distribution.

\medskip

\begin{theorem}\label{te1} Let $1<p<\infty$, $K\subset \mathrm{int}\mathbb{R}^n_+$ compact, and $L=\sum_{k=0}^m a_k(I-\Delta_a)^k$, $a_k\in \mathbb{R}$, $k=1, \dots , m$ be defined in a bounded open neighborhood $O$ of $K$ such that $\overline{O}\subset \mathrm{int}\mathbb{R}^n_+$ . Let $u\in L^p_a(O\setminus K)$ is a solution to $Lu=0$ in $O\setminus K$. If $N_{p',a,2m}(K)=0$, $u$ can be extended to $\tilde{u}\in L^p_a(O)$ such that $L\tilde{u}=0$ in $O$ in weak sense.
\end{theorem}

\medskip

\begin{remark}\label{r2}
{\rm (1) Let $1<p<\infty$, $\nu>0$. If $N_{p,a,\nu}(K)=0$, then $\lambda(K)=0$.\\
Indeed, let $\varepsilon>0$ be arbitrary and $f=G_{a,\nu}*_a g$ such that $\|g\|_{p',a}\le \varepsilon$, $f\equiv 1$ on $K$. By \eqref{gn}
$$\varepsilon^p\ge \|G_{a,\nu}*_a g\|_{p',a}^p\ge \int_K |G_{a,\nu}*_a g|^pd\mu_a=\mu_a(K).$$ Thus $\lambda_a(K)=0$ and so $\lambda(K)=0$.

\noindent (2) If $f\in C_0^\infty(\mathrm{int}\mathbb{R}^n_+)$ and $g\in C_0^\infty(\mathrm{int}\mathbb{R}^n_+)'$, that is in the dual space, then
\begin{equation}\label{adj} \langle Lf,g\rangle_a=\langle f,Lg\rangle_a.\end{equation}
In one dimension for $L=B_\alpha$ it is \cite[(2.4)]{p}. Similarly to this case if $f \in C_0^\infty(\mathrm{int}\mathbb{R}^n_+)$ and $g$ is smooth enough, then integration by parts implies the result. For general elements of  the dual space we extend  $Lg$ by the formula above, that is $Lg f:=\langle Lf,g\rangle_a$.
}\end{remark}

\medskip

\proof (of Theorem \ref{te1})
Let $O$ be a bounded open neighborhood of $K$ in $\mathbb{R}^n_+$ and $u\in L^p_a(O\setminus K)$ for which $Lu=0$ in $O\setminus K$. Since $N_{p',a,2m}(K)=0$, for all $\varepsilon>0$ there is a $\varphi=G_{a,2m}*_a f \in \mathcal{S}$ such that $\varphi\equiv 1$ on a neighborhhood $U\subset O$ of $K$ and $\|f\|_{p',a}\le \varepsilon$. Let $g\in C_0^\infty(O)$. Then $(1-\varphi)g \in C_0^\infty(O\setminus K)$. In view of (1) of Remark \ref{r2} $u$ is a.e. defined, so it can be handled as a distribution. Since $(1-\varphi)g \in C_0^\infty(O)$, by the assumption
$$\langle u,L(1-\varphi)g\rangle_a=0.$$
This implies that
$$\left|\langle u,Lg\rangle_a\right|=\left|\langle u,L\varphi g\rangle_a\right|\le \|u\|_{p,a,O}\|L\varphi g\|_{p',a}\le c\|u\|_{p,a,O}\|\varphi g\|_{W^m_{p',a}},$$
where $c=c(L)$ depends on the coefficients of $L$.\\
Applying Lemma \ref{l1} and the inversion of Bessel potential (for formulae see e.g. \cite[Theorem 1]{dls}) we have
$$\left|\langle u,Lg\rangle_a\right|\le c\|u\|_{p,a,O}\|\varphi g\|_{p',a,\nu}\le c \|\varphi g\|_{p',a}\le c\|g\|_\infty \|f\|_{p',a},$$
where $c$ depends on $u$ and $L$. Since $\varepsilon$ was arbitrary, for all $g\in C_0^\infty (O)$   $\langle u,Lg\rangle_a=0$, so in view of \eqref{adj} $u$ is a weak solution on $O$.

\medskip

The fundamental solution for the Laplace-Bessel operator, that is
$$\Delta_a E=\delta_a, \ws \ws\mbox{where}  \ws \langle \delta_a, \varphi\rangle_a=\varphi(0), \ws \ws \varphi \in \mathcal{S}_e,$$
is
\begin{equation}\label{E}E(x)=\left\{\begin{array}{ll}c(n,a)\ln x , n+|a|=2\\ c(n,a)|x|^{2-n-|a|}, n+|a|>2, \end{array}\right.\end{equation}
see e.g. \cite[Theorem 93 (page 324)]{ss}.

\begin{cor}With the notation above, let $L=\Delta_a$, that is $\Delta_au=0$ on $O\setminus K$, and $u\in L^p_a(O\setminus K)$. Let $2\le \nu< \frac{n+|a|}{p'}$. Then $u$ can be extended to $\tilde{u}\in L^p_a(O)$ such that $L\tilde{u}=0$ in $O$ in weak sense if and only if $C_{p',a,\nu}(K)=0$.
\end{cor}

\proof If $C_{p',a,\nu}(K)>0$, there is a nonzero measure $\mu\in\mathcal{M}(K)$ such that $G_{a.\nu}*_a\mu \in L^p_a(\mathbb{R}^n_+)$. In wiev of \eqref{E}, $E*_a\mu \in L^p_{a,loc}(\mathbb{R}^n_+)$. Since $E*_a\mu$ is a solution to $\Delta_a u=0$ in $O\setminus K$, $K$ is not removable.

On the other hand we prove that if $C_{p',a,\nu}(K)=0$, then $N_{p',a,2m}(K)=0$ too. According to Definition \ref{cap} for an $\varepsilon>0$ there is a nonnegative function $f\in L^p_a(\mathbb{R}^n_+)$ such that $G_{a,\nu}*_a f\ge 1$ on a neighborhood of $K$ and $\|f\|_{p,a}^p\le \varepsilon$.\\
Define a function $h\in C^\infty(\mathbb{R}_+)$, $0\le h\le 1$ and $h(t)=0$ on $t\in\left[0,\frac{1}{2}\right]$; $h(t)=1$ if $t\ge 1$. Taking into consideration that $G_{a,\nu}*_a f\ge 0$, $\varepsilon\ge\int_{\mathbb{R}_+}(G_{a,\nu}*_a f)^px^adx\ge \int_{G_{a,\nu}*_a f\ge\frac{1}{2}}\frac{1}{2^p}x^adx$. Thus $\int_{\mathbb{R}_+}(h(G_{a,\nu}*_a f))^px^adx\le \int_{G_{a,\nu}*_a f\ge\frac{1}{2}}x^adx\le 2^p\varepsilon$.\\
Noticing that $h(G_{a,\nu}*_a f)$ fulfils the requirements of Definition \ref{N}, since $\varepsilon$ is arbitrary, the statement is proved.

\section{Maximal measure and a Wolff type inequality}

Bessel maximal function was introduced and examined e.g. in cf. e.g. \cite{g}, see also the references therein. The boundedness of the maximal operator in some Morrey spaces is studied and applied to prove a Hardy-Littlewood-Sobolev type theorem in \cite{gu2}. The maximal measure presented below has proved useful in formulating a Wolff type inequality which is the main tool of the next section. Wolff type inequalities can be applied in different situations, for instance in martingale theory, see \cite{aclm} or deducing trace inequalities or characterize the trace measures via Wolff's inequality, see e.g. \cite{v}, \cite{v1}, \cite{t} and the references therein.

Below we define the maximal measure with respect to Bessel convolution.

\medskip

\begin{defi}$$M_a\mu(x):=\sup_{r>0}\frac{1}{\lambda_a(B_+(0,r))}\chi_{B_+(0,r)}*_a\mu(x).$$

\end{defi}

\medskip

Since
$$\lambda_a(B_+(0,r))=cr^{n+|a|},$$
we define the fractional maximal measure as
\begin{equation}\label{M}M_{a,d}\mu(x):= \sup_{r>0}\frac{\chi_{B_+(0,r)}*_a\mu(x)}{r^{n+|a|-d}},\end{equation}
and the truncated one as
\begin{equation}\label{Md}M_{a,d,b}\mu(x):=\sup_{0<r\le b}\frac{\chi_{B_+(0,r)}*_a\mu(x)}{r^{n+|a|-d}}.\end{equation}

\medskip

\begin{lemma}\label{lIm} $$I_\beta*_a\mu(x)=c\int_0^\infty \frac{\chi_{B_+(0,r)}*_a\mu(x)}{r^{n-\beta}}\frac{dr}{r}.$$
$$I_\beta\chi_{B_+(0,\delta)}*_a\mu(x)=c\int_0^\delta \frac{\chi_{B_+(0,r)}*_a\mu(x)}{r^{n-\beta}}\frac{dr}{r}+c\frac{\chi_{B_+(0,\delta)}*_a\mu(x)}{\delta^{n-\beta}}.$$
\end{lemma}

\proof
Let $d\Theta_{x,a}(z)=T^x\mu(z)z^{a}dz$. Changing the order of integration we get
$$\int_0^\delta \frac{\chi_{B_+(0,r)}*_a\mu(x)}{r^{n-\beta}}\frac{dr}{r}=\int_0^\delta \frac{1}{r^{n-\beta+1}}\int_{B_+(0,r)}1d\Theta_{x,a}(z)dr$$ $$=\int_{B_+(0,\delta)} \int_{|z|}^\delta \frac{1}{r^{n-\beta+1}}drd\Theta_{x,a}(z)=c\int_{B_+(0,\delta)}I_\beta(z)d\Theta_{x,a}(z)-c\int_{B_+(0,\delta)}\frac{1}{\delta^{n-\beta}}d\Theta_{x,a}(z)$$ $$=cI_\beta\chi_{B_+(0,\delta)}*_a\mu(x)-c\frac{\chi_{B_+(0,\delta)}*_a\mu(x)}{\delta^{n-\beta}}.$$

\medskip

\note
$$I_\beta^\delta*_a\mu(x):=c\int_0^\delta \frac{\chi_{B_+(0,r)}*_a\mu(x)}{r^{n-\beta}}\frac{dr}{r}.$$

\medskip

\begin{theorem}\label{Th2} Let $1\le p <\infty$, $0<\nu<n+|a|$, $\delta>0$ arbitrary. Then for all positive measure $\mu$, there are constants $c=c(n,a,\delta)$ (are not the same at each occurrence) such that
\begin{equation}\label{IM1} \|I_{\nu-|a|}*_a\mu\|_{p,a}\le c \|M_{a,\nu}\mu\|_{p,a}\end{equation}
and
\begin{equation}\label{IM2}\|I_{\nu-|a|}^\delta*_a\mu\|_{p,a}\le c \|M_{a,\nu,\delta}\mu\|_{p,a}.\end{equation}

\end{theorem}

\medskip

\note
$$H_s^\mu:=\{x : I_{\nu-|a|}*_a\mu(x)>s\}, \ws \ws K_s^\mu:=\{x: M_{a,\nu}\mu(x)> s\}.$$
$$^1H_s^\mu:=\{x : I_{\nu-|a|}^1*_a\mu(x)>s\}, \ws \ws ^1K_s^\mu:=\{x: M_{a,\nu,1}\mu(x)> s\}.$$

\medskip

\begin{lemma} There is a $\varrho>1$ and a $b>0$, such that for all $s>0$ and $\varepsilon \in (0,1]$,
\begin{equation}\label{Hs}\lambda_a(H_{\varrho s}^\mu)\le b\varepsilon^{\frac{n+|a|}{n+|a|-\nu}}\lambda_a(H_{s}^\mu)+\lambda_a(K_{\varepsilon s}^{\mu}).\end{equation}
Similarly
\begin{equation}\label{Hs1}\lambda_a(^1H_{\varrho s}^\mu)\le b\varepsilon^{\frac{n+|a|}{n+|a|-\nu}}\lambda_a(^1H_{s}^\mu)+\lambda_a(^1K_{\varepsilon s}^{\mu}).\end{equation}
\end{lemma}

\proof  By lower semicontinuity we can take Whithney's decomposition of $H_s^\mu$, i.e. $H_s^\mu=\cup_{i=1}^\infty Q_i$, where $Q_i$-s are dyadic cubes, $\mathrm{int}Q_i\cap\mathrm{int}Q_j=\emptyset$ if $i\neq j$ and $\mathrm{diam}Q_i\le \mathrm{dist}(Q_i, (H_s^\mu)^c)\le 4\mathrm{diam}Q_i$. (Dyadic cubes means cubes with side $2^{-k}$, $k\in\mathbb{Z}$, whose vertices belong to the lattice $\{m2^{-k}: m\in\mathbb{Z}^n\}$. For Whithney's decomposition see \cite[page 16, Theorem 3]{st}.) In addition, to prove \eqref{Hs1} if $\mathrm{diam}Q_i\ge \frac{1}{8}$, then we decompose it to subcubes with diameter is between $\frac{1}{16}$ and $\frac{1}{8}$, and we consider this new sequence of cubes.

Let $Q$ be an element of this decomposition. Let $x\in Q$ be arbitrary, denote the center of $Q$ by $x_c$ and let  $d:=\mathrm{diam}Q$. Let $G:=B(x_c, 6d)$, $B=B(x, 8d)$, that is $Q\subset G\subset B$. Let $\mu=\mu_1+\mu_2$, where $\mu_1=\mu|_G$.

At first we deal with  $I_{\nu-|a|}*_a\mu_2$. To this we estimate $T^t\chi_{B_+(0,r)}(x)$. We can assume, that $r>\frac{11}{2}$, otherwise $\mathrm{supp}\mu_2\cap B_+(x,r)=\emptyset$. Let $x_1\in (H_s^\mu)^c$ such that $\mathrm{dist}(x_1,Q)<4d$. Then by \eqref{kf} if $t\in B_+(x,r)$,
$$T^t\chi_{B_+(0,r)}(x)\le  c \prod_{i=1}^n\min\left\{1,\left(\frac{r}{x_{i}}\right)^{a_i}\right\}\le  c \prod_{i=1}^n\min\left\{1,\left(\frac{r}{x_{1,i}}\right)^{a_i}\right\}=:P.$$
Indeed, it is obvious if $x_{1,i}\le x_{i}$. If $x_{1,i}, x_{i}<r$, then the minimum is 1 in both cases. If $x_{i}< r \le x_{1,i}$,  then $\min\left\{1,\frac{r}{x_{i}}\right\}=1$  thus we have to show that $1\le c \frac{r}{x_{1,i}}$. Since $x_{1,i}<x_i+5d$, $\frac{r}{x_{1,i}}>\frac{r}{x_{i}+5d}>\frac{r}{r+5d}\ge \frac{\frac{11}{2}d}{\frac{11}{2}d+5d}= \frac{11}{21}$.\\
If $\frac{11}{2}d< r <x_{i}<  x_{1,i}$, then $\frac{r}{x_{1,i}}>\frac{r}{x_{i}+\frac{20}{11}r}>\frac{11}{21}\frac{r}{x_{i}}$.

In view of \eqref{kf1} if $t\in T_+(x_1,2r)$
$$P\le cT^t\chi_{[0,4r)^n}(x_1).$$
Thus
$$\chi_{B_+(x,r)}*_a\mu_2\le c\chi_{[0,4r)^n}*_a\mu_2(x_1)\le c \chi_{B_+(0, 4\sqrt{n}r)}*_a\mu_2(x_1).$$
Recalling that $x_1\in (H_s^\mu)^c$, we have
$$I_{\nu-|a|}*_a\mu_2(x)\le c\int_{\frac{11}{2}d}^\infty \frac{\chi_{B_+(0, 4\sqrt{n}r)}*_a\mu_2(x_1)}{(4\sqrt{n}r)^{n+|a|-\nu}}dr=c I_{\nu-|a|}*_a\mu_2(x_1)\le c s.$$

Now we choose $\varrho$ so that $I_{\nu-|a|}*_a\mu_2(x)\le \frac{\varrho s}{2}$, which implies that
\begin{equation}\label{mu1}H_{\varrho s}^\mu\cap Q\subset H_{\frac{\varrho s}{2}}^{\mu_1}\cap Q.\end{equation}
If the diameter of $Q$ was originally less then $\frac{1}{8}$, then the whole construction is contained in a ball of radius less than one, so the same chain of ideas leads to
\begin{equation}\label{mu11}^1H_{\varrho s}^\mu\cap Q\subset \ws ^1H_{\frac{\varrho s}{2}}^{\mu_1}\cap Q.\end{equation}
If there is no $x_1\in (^1H_s^\mu)^c$ such that $\mathrm{dist}(Q, x_1)\le 4d$, then $\mathrm{diam}Q>\frac{1}{16}$. Let $x_0\in Q\cap (^1K_{\varepsilon s})^c$. Then, recalling that $r>\frac{11}{2}d$, we have
$$I_{\nu-|a|}^1*_a\mu_2(x_0)=c\int_{\frac{11}{32}}^1\frac{\chi_{B_+(0,r)}*_a\mu_2(x_0)}{r^{n+|a|-\nu+1}}dr $$ $$\le \frac{32}{11}\int_0^1\frac{\chi_{B_+(0,r)}*_a\mu_2(x_0)}{r^{n+|a|-\nu}}dr \le c M_{a,\nu,1}(x_0)\le c \varepsilon s.$$
Thus if $Q\cap (^1K_{\varepsilon s})^c\neq\emptyset$, we can choose $\varrho$ again so that \eqref{mu11} is satisfied.

Let $x_0\in Q \cap (K_{\varepsilon s})^c$ again. According to \cite[Theorem 2 (c)]{g}
$$\lambda_a\left(Q\cap H_{\frac{\varrho s}{2}}^{\mu_1}\right)\le c \left(\frac{1}{\varrho s}\int_{\mathbb{R}^n_+}1d\mu_{1,a}(t)\right)^{\frac{n+|a|}{n+|a|-\nu}}=c\left(\frac{1}{\varrho s}\int_{G}1d\mu_{a}(t)\right)^{\frac{n+|a|}{n+|a|-\nu}}$$ $$\le c\left(\frac{1}{\varrho s}\int_{B}1d\mu_{a}(t)\right)^{\frac{n+|a|}{n+|a|-\nu}}=(*).$$
In view of Lemma \ref{tk} $B=\mathrm{supp}T^t\chi_{B_+(0,8d)(x_0)}\subset T_+(x_0,8d)$. Thus by \eqref{kf1}
$$(*)\le c\left(\frac{1}{\varrho s}\prod_{i=1}^n\max\left\{1, \left(\frac{x_{0,i}}{16d}\right)^{a_i}\right\}\int_{\mathbb{R}^n_+}T^t\chi_{[0,16d)^n}(x_0)d\mu_{a}(t)\right)^{\frac{n+|a|}{n+|a|-\nu}}$$
\begin{equation}\label{M1}\le c\left(\frac{M_{a,\nu}\mu(x_0)}{\varrho s}\right)^{\frac{n+|a|}{n+|a|-\nu}}(\prod_{i=1}^n\max\{1, \left(\frac{x_{0,i}}{16d}\right)^{a_i}\}))^{\frac{n+|a|}{n+|a|-\nu}}d^{n+|a|}.\end{equation}
Taking into consideration that $\lambda_a(Q)\sim \prod_{i=1}^n x_{c,i}^{a_i}\left(\frac{d}{\sqrt{n}}\right)^n\sim \prod_{i=1}^n\left(\frac{x_{c,i}}{d}\right)^{a_i}d^{n+|a|}$ and recalling that $|x_{0,i}-x_{c,i}|<d$, we have
$$\left(\prod_{i=1}^n\max\{1, \left(\frac{x_{0,i}}{16d}\right)^{a_i}\}\right)^{\frac{n+|a|}{n+|a|-\nu}}d^{n+|a|}\le c |Q|_a \left(\prod_{i=1}^n\left(\frac{x_{c,i}}{d}\right)^{a_i}\right)^{\frac{\nu}{n+|a|-\nu}}\le  c \lambda_a(Q).$$
Finally as $x_0\in (K_{\varepsilon s})^c$,
\begin{equation}\label{QH}\lambda_a\left(Q\cap H_{\frac{\varrho s}{2}}^{\mu_1}\right)\le b \varepsilon^{\frac{n+|a|}{n+|a|-\nu}}\lambda_a(Q).\end{equation}
Since $\mathrm{diam}Q\le \frac{1}{8}$, similarly we have
\begin{equation}\label{QH1}\lambda_a\left(Q\cap \ws ^1H_{\frac{\varrho s}{2}}^{\mu_1}\right)\le b \varepsilon^{\frac{n+|a|}{n+|a|-\nu}}\lambda_a(Q),\end{equation}
cf. \eqref{M1}.

That is if $Q\cap (K_{\varepsilon s})^c\neq\emptyset$ or $Q\cap (^1K_{\varepsilon s})^c\neq\emptyset$, then \eqref{QH} or \eqref{QH1}, respectively, is fulfilled, otherwise $Q\subset K_{\varepsilon s}\mu$.\\
Recalling \eqref{mu1} or \eqref{mu11} and adding over all $Q\in\{Q_i\}$, we obtain the required result.

\medskip

\proof (of Theorem \ref{Th2}) \eqref{IM1} and \eqref{IM2} can be derived by the same way from \eqref{Hs} and \eqref{Hs1}. Let us see the second one, say. Let $\delta=1$. Integrating \eqref{Hs1} and changing the variables we have
$$\frac{1}{\varrho^p}\int_0^{\varrho R}\lambda_a(^1H_{u}^\mu) u^{p-1}du\le b\varepsilon^{\frac{n+|a|}{n+|a|-\nu}}\int_0^R\lambda_a(^1H_{s}^\mu)u^{p-1}du+\frac{1}{\varepsilon^p}\int_0^{\varepsilon R}\lambda_a(^1K_{u}^{\mu})u^{p-1}du.$$
Supposing that $\mathrm{supp}\mu$ is bounded, all the integrals are finite. We choose $\varepsilon$ small enough so that $b\varepsilon^{\frac{n+|a|}{n+|a|-\nu}}\le \frac{1}{2\varrho^p}$. Then we have
$$\frac{1}{\varrho^p}\int_0^{\varrho R}\lambda_a(^1H_{u}^\mu) u^{p-1}du\le \frac{2}{\varepsilon^p}\int_0^{\varepsilon R}\lambda_a(^1K_{u}^{\mu})u^{p-1}du.$$
Letting $R \to \infty$,
$$\|I_{\nu-|a|}^1*_a\mu(x)\|_{a,p}\le\frac{\varrho}{\varepsilon}\|M_{a,\nu,1}\mu(x)\|_{a,p},$$
cf. e.g. \cite[(1.46)]{ss}. If $\mathrm{supp}\mu$ is not compact, then let $\mu_n:=\mu|_{B_+(0,n)}$. Since $\|M_{a,\nu,1}\mu_n\|_{a,p}\le \|M_{a,\nu,1}\mu\|_{a,p}$, we can apply the monotone convergence theorem.

\medskip

\begin{cor}\label{cor1} With the assumptions of Theorem \ref{Th2} we have
\begin{equation}\|I_{\nu-|a|}^\delta*_a\mu\|_{p,a}\le c \|G_{a,\nu}*_a\mu\|_{p,a}\le c  \|M_{a,\nu,\delta}\mu\|_{p,a}.\end{equation}
\end{cor}

\proof
\eqref{k1} and \eqref{k2} ensure the first inequality and
$$\|G_{a,\nu}*_a\mu\|_{p,a}\le \|I_{\nu-|a|}\chi_{B_+(0,\delta)}*_a\mu\|_{p,a}+\|e^{-\frac{|\cdot|}{2}}*_a\mu\|_{p,a}$$ $$\le \|I_{\nu-|a|}^\delta*_a\mu\|_{p,a}+\|M_{a,\nu,\delta}\mu\|_{p,a}+\|e^{-\frac{|\cdot|}{2}}*_a\mu\|_{p,a}.$$
Observe that
$$\frac{e^{-\frac{|\cdot|}{2}}*_a\chi_{B_+(0,r)}(x)}{r^{n+|a|}}$$ $$=\frac{1}{r^{n+|a|}}\int_{B_+(0,r)}\int_{[0,\pi)^n}e^{-\frac{1}{2}\sqrt{\sum_{i=1}^nx_i^2+t_i^2-2x_it_i\cos\vartheta_i}}d\sigma^{a}(\vartheta)d\lambda_a(t)$$ $$\ge \frac{1}{r^{n+|a|}}\int_{B_+(0,r)}e^{-\frac{1}{2}|x+t|}d\lambda_a(t)\ge c e^{-\frac{r}{2}}e^{-\frac{|x|}{2}}.$$
Thus
$$e^{-\frac{|x|}{2}}\le ce^{\frac{r}{2}}\frac{e^{-\frac{|\cdot|}{2}}*_a\chi_{B_+(0,r)}(x)}{r^{n+|a|}},$$
and so if $r\le \delta$,
$$e^{-\frac{|\cdot|}{2}}*_a\mu(x)\le c(r)e^{-\frac{|\cdot|}{2}}*_a\frac{\chi_{B_+(0,r)}(x)}{r^{n+|a|}}*_a\mu\le c(r,\nu)e^{-\frac{|\cdot|}{2}}*_aM_{a,\nu, \delta}\mu.$$
According to \eqref{Y}
$$\|e^{-\frac{|\cdot|}{2}}*_a\mu\|_{p,a}\le\|e^{-\frac{|\cdot|}{2}}\|_{1,a}\|M_{a,\nu,\delta}\mu\|_{p,a},$$
which, together with Theorem \ref{Th2}, implies the statement.

\medskip

\note
Denote by
$$b_k(x):=2^{k(n+|a|-\nu)}\chi_{B_+(0,2^{-k})}*_a\mu(x)$$
and by
$$c_k(x):=2^{k(n+|a|-p\nu)}\chi_{B_+(0,2^{-k})}*_a\mu(x).$$
The corresponding Wolff-function is
$$W_{a,\nu,p}^\mu(x):=\int_0^1\left(\frac{\chi_{B_+(0,r)}*_a\mu(x)}{r^{n+|a|-p\nu}}\right)^{p'-1}\frac{dr}{r}.$$

\medskip

\begin{remark}
{\rm In view of Lemma \ref{lIm} and \eqref{Md} we can observe that
\begin{equation}\label{b}I_{\nu-|a|}*_a\mu(x)\sim \|\{b_k(x)\}_{-\infty}^\infty\|_{l^1}, \ws \ws I_{\nu-|a|}^1*_a\mu(x) \sim \|\{b_k(x)\}_{0}^\infty\|_{l^1},\end{equation}
and
\begin{equation}\label{b2}M_{a,\nu, 1}\mu(x)\sim \|\{b_k(x)\}_{0}^\infty\|_{l^\infty}.\end{equation}}\end{remark}

\medskip

We are in position to prove the next Wolff type inequality.

\begin{theorem}\label{cw} $1< p , q < \infty$
$$\|G_{a,\nu}*_a\mu\|_{p',a} \sim \|\|\{b_k(x)\}_{0}^\infty\|_{l^q}\|_{p',a}$$ $$\sim \left(\int_{\mathbb{R}^n_+}\|\{c_k(x)\}_{0}^\infty\|_{l^{p'-1}}^{p'-1}d\mu_a(x)\right)^{\frac{1}{p'}}\sim \left(\int_{\mathbb{R}^n_+}W_{a,\nu,p}^\mu(x)d\mu_a(x)\right)^{\frac{1}{p'}}.$$
\end{theorem}

\proof

Corollary \ref{cor1}, \eqref{b} and \eqref{b2} ensure that
$$\|G_{a,\nu}*_a\mu\|_{p,a}\le c\|M_{a,\nu,1}\mu\|_{p,a}\le  c\|\|\{b_k(x)\}_{0}^\infty\|_{l^\infty}\|_{p,a}$$ $$ \le\|\|\{b_k(x)\}_{0}^\infty\|_{l^p}\|_{p,a}\le\|\|\{b_k(x)\}_{0}^\infty\|_{l^1}\|_{p,a}\le c\|I_{\nu-|a|}^2*_a\mu\|_{p,a}\le c \|G_{a,\nu}*_a\mu\|_{p,a}.$$
To prove the Wolff type inequality we have
$$\|\|\{b_k(x)\}_{0}^\infty\|_{l^{p'}}\|_{p',a}^{p'}=\int_{\mathbb{R}^n_+}\sum_{k=0}^\infty b_k(t)^{p'}d\lambda_a(t)$$ $$=\sum_{k=0}^\infty 2^{k(n+|a|-\nu)p'}\int_{\mathbb{R}^n_+}\chi_{B_+(0,2^{-k})}*_a\mu(t)(\chi_{B_+(0,2^{-k})}*_a\mu)^{p'-1}(t)d\lambda_a(t)$$ $$=\sum_{k=0}^\infty 2^{k(n+|a|-\nu)p'}I_k.$$
$$I_k=\int_{\mathbb{R}^n_+}\int_{\mathbb{R}^n_+}T^x\chi_{B_+(0,2^{-k})}(t)d\mu_a(x)(\chi_{B_+(0,2^{-k})}*_a\mu)^{p'-1}(t)d\lambda_a(t)$$ $$=\int_{\mathbb{R}^n_+}\int_{\mathbb{R}^n_+}\int_{\mathbb{R}^n_+}K(x,t,z)\chi_{B_+(0,2^{-k})}(z))d\lambda_a(z)d\mu_a(x)(\chi_{B_+(0,2^{-k})}*_a\mu)^{p'-1}(t)d\lambda_a(t)$$ $$=\int_{\mathbb{R}^n_+}\int_{B_+(0,\frac{1}{2^k})}\int_{\mathbb{R}^n_+}K(x,t,z)(\chi_{B_+(0,2^{-k})}*_a\mu)^{p'-1}(t)d\lambda_a(t)d\lambda_a(z))d\mu_a(x).$$
Since $0<z_i<2^{-k}$, recalling that $t_i\in (|x_i-z_i|,x_i+z_i)$,\\ $(\chi_{B_+(0,2^{-k})}*_a\mu)(t)\sim (\chi_{B_+(0,2^{-k})}*_a\mu)(x)$. As $K$ is a reproducing kernel
$$I_k\sim \lambda_a(B_+(0,2^{-k}))\int_{\mathbb{R}^n_+}(\chi_{B_+(0,2^{-k})}*_a\mu)^{p'-1}(x)d\mu_a(x).$$
Thus recalling that $\lambda_a(B_+(0,2^{-k})) \sim 2^{-k(n+|a|)}$ by Fubini's theorem again
$$\|\|\{b_k(x)\}_{0}^\infty\|_{l^q}\|_{p',a}^{p'}\sim \int_{\mathbb{R}^n_+}\sum_{k=0}^\infty c_k(x)^{p'-1}d\mu_a(x).$$

\medskip

\begin{remark}\label{nem1} {\rm In view of Corollary \ref{cor1}, instead of $M_{a,\nu,1}\mu$ we can consider $M_{a,\nu,\delta}\mu$ with the corresponding sequence $\{b_k^\delta(x)\}$.}
\end{remark}

\section{Metric properties}

Applying the previous section, below we investigate some "metric" properties of B-$p$ capacity. Since the Bessel-translation is not a geometric congruence, we need a special "Lipschitz"- condition. It is also necessary to introduce the notion of "weighted Hausdorff measure", to examine Cantor-type sets.

At the beginning of this section let us recall that for $1<p<\infty$ the B-$p$ capacity of $K\subset \mathbb{R}^n_+$ is
\begin{equation}\label{dc}C_{a,\nu,p}^{\frac{1}{p}}(K):=\sup\{\mu_a(K) : \mu \in \mathcal{M}(K), \ws \|G_{a,\nu}*_a\mu\|_{p,a}\le 1\},\end{equation}
and the B-$p$ capacity is non trivial if $1<p< \frac{n+|a|}{\nu}$.

\medskip

\subsection{B-Lipschitz mappings}
The next Lipschitz type property is corresponding to the Bessel translation.

\medskip

\begin{defi}Let  $\Phi:\mathbb{R}^n_+ \to \mathbb{R}^n_+$. Consider $z(x,t,\vartheta)=z_1(x,t,\vartheta), \dots, z_n(x,t,\vartheta)$, where $z_k(x,t,\vartheta)=x_k-t_k\cos\vartheta_k+it_k\sin\vartheta_k$. Let $\Psi:\mathbb{C}^n_+ \to \mathbb{C}^n_+$ such that $\Psi(z)_k(x,t\vartheta)=\Phi(x)_k-\Phi(t)_k\cos\vartheta_k+i\Phi(t)_k\sin\vartheta_k$. We say that $\Phi$ fulfils the B-Lipschitz condition with B-Lipschitz constant $L$ if for a.e. $\vartheta \in [0\pi)^n$
\begin{equation}\label{BL}|\Psi(z)(\vartheta)|\le L |z(\vartheta)|.\end{equation}

\end{defi}

\medskip

\begin{remark} {\rm Of course, linear mappings posses the B-Lipschitz property.\\
For  $\vartheta=0$ \eqref{BL} gives back the standard Lipschitz condition, and for $\vartheta_k=\pi$ $k=1, \dots , n$ \eqref{BL} means that $|\Phi(x)+\Phi(t)|\le L|x+t|$.}
\end{remark}

\medskip

\noindent {\bf Example.} Let $f: \mathbb{R}_+ \to \mathbb{R}_+$ be a Lipschitz function. Let $K\subset \mathrm{int}\mathbb{R}_+^n$ compact and $\Phi: K\to \mathbb{R}_+^n$; $\Phi(x)=f(|x|)x$. On $K$ $\Phi$ also fulfils the Lipschitz condition with constant $\tilde{L}(K)$ and $\frac{\Phi(x)_k}{x_k}\le M(K)$. Let $$G(\vartheta):=L^2 L |z(\vartheta)|^2-|\Phi(z)(\vartheta)|^2$$ $$=\sum_{k=1}^n L^2(x_k^2+t_k^2)-\Phi(x)_k^2-\Phi(t)_k^2-2\sum_{k=1}^n \cos \vartheta_k(L^2x_kt_k-\Phi(x)_k\Phi(t)_k).$$
If $L>M(K)$, we have
$$\frac{\partial G}{\partial \vartheta_k}=0 \ws \ws \mbox{iff} \ws \vartheta_k=0 \ws \mbox{or} \ws \vartheta_k=\pi.$$
It can be readily seen, that the Hessian is positive (negative) definite if $\vartheta=0$ ($\vartheta_k=\pi$ $k=1, \dots , n$), respectively.  Thus the Lipschitz property of $\Phi(x)$ implies that the B-Lipschitz condition fulfils for all $\vartheta \in [0\pi)^n$.

\medskip

\begin{theorem} Let $\nu>0$, $1<p< \frac{n+|a|}{\nu}$. Let $E\subset \mathbb{R}^n_+$ and $\Phi: E\to \mathbb{R}^n_+$ is a B-Lipschitz mapping with B-Lipschitz constant $L$. Then
$$C_{p,a,\nu}(\Phi(E))\le c \ws C_{p,a,\nu}(E),$$
where $c$ depends only on $n, p, a, L$.
\end{theorem}

\proof By standard arguments it is enough to prove for any $K\subset \mathbb{R}^n_+$, compact. Let $\mu \in \mathcal{M}(\Phi(K))$ Then by \cite[Lemma 5.2.2]{ah} there is a $\mu_\Phi \in \mathcal{M}(K)$ such that
$$\int_{\Phi(K)}W_{a,\nu,p}^\mu(y)d\mu_a(y)=\int_K W_{a,\nu,p}^\mu(\Phi(x))\Phi^a(x)d\mu_\Phi(x)$$ $$=\int_K \int_0^1\left(\frac{\int_K T^{\Phi(u)}\chi_{B_+(0,r)}(\Phi(x))\Phi^a(u)d\mu_\Phi(u)}{r^{n+|a|-\nu p}}\right)^{p'-1}\frac{dr}{r}\Phi^a(x)d\mu_\Phi(x)$$
$$\ge \int_K \int_0^1\left(\frac{\chi_{B_+\left(0,\frac{r}{L}\right)}*_a\mu_\Phi(x)}{r^{n+|a|-\nu p}}\right)^{p'-1}\frac{dr}{r}\Phi^a(x)d\mu_\Phi(x)$$ $$=c\int_K \int_0^{\frac{1}{L}}\left(\frac{\chi_{B_+\left(0,r\right)}*_a\mu_\Phi(x)}{r^{n+|a|-\nu p}}\right)^{p'-1}\frac{dr}{r}\Phi^a(x)d\mu_\Phi(x),$$
where $c=c(n,p,a,L)$.\\
This implies that
$$\int_{\Phi(K)}W_{a,\nu,p}^\mu(y)d\mu_a(y)>c\int_{K}W_{a,\nu,p}^{\mu_\Phi}(x)\Phi^a(x)d\mu_\Phi(x).$$
Indeed, if $L\le 1$, we have immediately the inequality above, if $L>1$ we have to consider Remark \ref{nem1} with $\delta=\frac{1}{L}$, which leads again to the inequality above. According to Theorem \ref{cw} it proves the statement, cf. the definition above.\\

\medskip

\subsection{Coverings}

In the next subsections coverings in Bessel-weighted space are introduced. Since Bessel-convolution lives in a weighted space, B-$p$ capacity of a set depends also on the location of the set. As capacity is in close connection with Hausdorff measure, in the next subsection we extend this notion to weighted spaces as well.

\medskip

\note Let $K\subset \mathbb{R}^n_+$ compact.

(1) Let $\mathcal{A}(r)$ be the minimal number of balls of radius $r$ required to cover $K$.

(2) $$m^a:=m^a(u,r,K)=\max\{x^a : x\in \overline{B(u,r)}\cap K\}.$$
$$\mathcal{B}(r):=\inf\{\sum_{j=1}^{\mathcal{A}(r)}m_j^a : K\subset\cup_{j=1}^{\mathcal{A}(r)}B(u_j,r) \ws \mbox{is a minimal covering}\}.$$

\begin{remark}\label{r5} {\rm

(1) $\mathcal{A}(r)$ is obviously decreasing. \\
(2) Let $\cup_{j=1}^{\mathcal{A}(r)}B(u_j,r)$ be a minimal covering of $K$. Then there is a constant $C_n$ such that any point of $K$ belongs to at most $C_n$ balls. Indeed, let $C_n=C_n(q)$ be the minimal number of balls of radius $q\le\frac{1}{2}$ which covers the unite ball. Suppose that there is a point $x\in K$ which belongs to $C_n(q)+1$ balls. Then $B\left(x, \frac{r}{q}\right)$ contains all these balls and it can be covered by $C_n(q)$ balls of radius $r$ which contradicts with minimality, cf. e.g. \cite[page 145]{ah}.\\
(3) Let $r_k:=\frac{1}{2^k}$, $r_{k+1}\le r\le r_k$ and $\frac{1}{4}\le q:=\frac{1}{2^{k+2}r}\le \frac{1}{2}$. Let $\{B(v_j,r_1)\}_{j=1}^{\mathcal{A}(r)}$ and $\{B(u_i,r_{k+2})\}_{i=1}^{\mathcal{A}(r_{k+2})}$ be minimal coverings of $K$ with the corresponding points $\{m_j\}_{j=1}^{\mathcal{A}(r)}$ and  $\{M_i\}_{i=1}^{\mathcal{A}(r_{k+2})}$, respectively. Any $m_j$ belongs to a ball, $B(u_i,r_{k+2})$, and so $m_j^a \le M_i^a$. Since at most $C_n\left(\frac{1}{q+2}\right)$ maximum points ($m_j$) can be in the same ball, $\sum_{j=1}^{\mathcal{A}(r)}m_j^a\le C_n\left(\frac{1}{q+2}\right)\sum_{j=1}^{\mathcal{A}(r_{k+2})}M_j^a$, thus $\mathcal{B}(r)\le C_n\left(\frac{1}{q+2}\right)\mathcal{B}(r_{k+2})$.  Repeating the chain of ideas with $r_{k-1}$ and $q_1=2^{k-1}r$, we have $\frac{1}{C(n,q_1)}\mathcal{B}(r_{k-1})\le \mathcal{B}(r)$.}
\end{remark}

\medskip

\begin{theorem}\label{tB} As above, let $\nu>0$, $1<p< \frac{n+|a|}{\nu}$. Then
$$C_{a,\nu,p}^{\frac{1}{p}}(K)\le c \left(\int_0^1\left(\mathcal{B}(r)r^{n-p\nu}\right)^{1-p'}\frac{dr}{r}\right)^{1-p},$$
where $c=c(n,a)$.
\end{theorem}

\proof
Let $\mu \in \mathcal{M}(K)$ and $r_k=\frac{1}{2^k}$, as above. According to Corollary \ref{cor1} and \eqref{b} we have
$$\|G_{a,\nu}*_a\mu\|_{p',a}^{p'}\ge \|I_{\nu-|a|}^1*_a\mu\|_{p',a}^{p'}$$ $$\ge c\int_{\mathbb{R}^n_+}\left(\sum_{k=0}^\infty 2^{k(n+|a|-\nu)}\chi_{B_+\left(0,r_k\right)}*_a\mu(x)\right)^{p'}d\lambda_a(x)$$
$$\ge c \sum_{k=0}^\infty \frac{1}{r_k^{p'(n+|a|-\nu)}}\int_{\mathbb{R}^n_+}\left(\chi_{B_+\left(0,r_k\right)}*_a\mu(x)\right)^{p'}d\lambda_a(x),$$
where the last inequality follows from the monotone convergence theorem.
Let $K\subset \cup_{j=1}^{\mathcal{A}(r_{k+1})}B(u_j,r_{k+1})$ be a minimal covering. Recalling that $\chi_{B_+\left(0,r_k\right)}*_a\mu(x)=\int_{\mathbb{R}^n_+}T^t\chi_{B_+\left(0,r_k\right)}(x)d\mu_a(t)$, if $t\in K$, $x\in U_k$, where $U_k$ is the neighbourhood of $K$ of radius $r_k$, otherwise $T^t\chi_{B_+\left(0,r_k\right)}(x)=0$. Noticing that $\cup_{j=1}^{\mathcal{A}(r_{k+1})}B(u_j,r_{k+1})\subset U_k$, in view of Remark \ref{r5} we have
$$\int_{\mathbb{R}^n_+}\left(\chi_{B_+\left(0,r_k\right)}*_a\mu(x)\right)^{p'}d\lambda_a(x)$$ $$\ge \frac{1}{C_n}\sum_{j=1}^{\mathcal{A}(r_{k+1})}\int_{B(u_j,r_{k+1})}\left(\int_{\mathbb{R}^n_+}T^t\chi_{B_+\left(0,r_k\right)}(x)d\mu_a(t)\right)^{p'}d\lambda_a(x)=:(**).$$
Considering again the support of $T^t\chi_{B_+\left(0,r_k\right)}$, by H\" older's inequality we have
$$(**)\ge \frac{1}{C_n}\sum_{j=1}^{\mathcal{A}(r_{k+1})}\lambda_a(B(u_j,r_{k+1}))^{-\frac{p'}{p}}$$ $$\times \left(\int_{B(u_j,r_{k+1})}\int_{B_+(x,r_k)}T^t\chi_{B_+(0,r_k)}(x)d\mu_a(t)d\lambda_a(x)\right)^{p'}=:(***).$$
Since $x\in B\left(u_j,r_{k+1}\right)$, $B(x,r_k)\supset B(u_j,r_{k+1})$. As $T^t\chi_{B_+\left(0,r_k\right)}(x)$ is continuous (actually it belongs to the $\mathrm{Lip}(\frac{1}{2})$ class) on  $B_+\left(x,r_k\right)$, $T^t\chi_{B_+\left(0,r_k\right)}(x)\ge c T^t\chi_{B\left(0,r_k\right)}(u_j)$. According to Lemma \ref{tk}, if $t\in T_+\left(u_j,\frac{r_k}{4\sqrt{n}}\right)=:T_{k,j}$, then $T^t\chi_{B\left(0,r_k\right)}(u_j)\ge c \prod_{i=1}^n \min\left\{1, \left(\frac{r_k}{ u_{j,i}}\right)^{a_i}\right\}$. Thus we have
$$ (***)$$ $$\ge  \frac{c}{C_n} \sum_{j=1}^{\mathcal{A}(r_{k+1})}\lambda_a(B\left(u_j,r_{k+1}\right))^{-\frac{p'}{p}}\mu_a^{p'}(B\left(u_j,r_{k+1}\right))$$ $$\times\left(\int_{T_{k,j}}\prod_{i=1}^n \min\left\{1, \left(\frac{r_k}{u_{j,i}}\right)^{a_i}\right\}d\lambda_a(x)\right)^{p'}.$$
Estimating the last integral we have
$$\int_{T_{k,j}}\prod_{i=1}^n \min\left\{1, \left(\frac{r_k}{u_{j,i}}\right)^{a_i}\right\}d\lambda_a(x)$$ $$\ge \prod_{i=1}^n\int_{u_{j,i}-\frac{r_{k+1}}{4\sqrt{n}}}^{u_{j,i}+\frac{r_{k+1}}{4\sqrt{n}}}\min\left\{1, \left(\frac{r_k}{u_{j,i}}\right)^{a_i}\right\}x_i^{a_i}dx_i=\prod_{i=1}^nI_{j,i}.$$
If $r_k>u_{j,i}$,
$$I_{j,i}\ge c\int_0^{\frac{r_{k+1}}{4\sqrt{n}}}x_i^{a_i}dx_i\ge c r_{k+1}^{a_i+1}.$$
If $r_k\le u_{j,i}$,
$$I_{j,i}\ge c \int_{u_{j,i}-\frac{r_{k+1}}{4\sqrt{n}}}^{u_{j,i}+\frac{r_{k+1}}{4\sqrt{n}}}r_{k+1}^{a_i}dx_i\ge cr_{k+1}^{a_i+1}.$$
So we have
$$(***)\ge c\sum_{j=1}^{\mathcal{A}(r_{k+1})}r_k^{n(1-p')}m_j^{a(1-p')}\mu_a^{p'}(B(u_j,r_{k+1}))r_k^{(n+|a|)p'}$$ $$=cr_k^{n+|a|p'}\sum_{j=1}^{\mathcal{A}(r_{k+1})}m_j^{a(1-p')}\mu_a^{p'}(B(u_j,r_{k+1})).$$
By H\"older's inequality
$$\mu_a(K)\le \left(\sum_{j=1}^{\mathcal{A}(r_{k+1})}\mu_a(B(u_j,r_{k+1}))^{p'}m_j^{a(1-p')}\right)^{\frac{1}{p'}}\left(\sum_{j=1}^{\mathcal{A}(r_{k+1})}m_j^{a}\right)^{\frac{1}{p}}.$$
Thus
$$\|G_{a,\nu}*_a\mu\|_{p',a}^{p'}\ge c\mu_a(K)^{p'}\sum_{k=0}^\infty \left(r_k^{n-\nu p}\mathcal{B}(r_{k+1})\right)^{1-p'}.$$
Taking into consideration Remark \ref{r5} we have
$$\|G_{a,\nu}*_a\mu\|_{p',a}^{p'}\ge c\mu_a(K)^{p'}\int_0^1\left(\mathcal{B}(r)r^{n-p\nu}\right)^{1-p'}\frac{dr}{r}.$$
Comparing it with \eqref{dc} the proof is finished.

\subsection{Hausdorff measure with Bessel external field}

\begin{defi}\label{Hau} Let $h$ be an increasing function on $\mathbb{R}_+$ with $h(0)=0$. Let $E\subset \mathbb{R}_+^n$, $\varrho>0$.
\begin{equation}\Lambda^{\varrho}_{h,a}(E):=\inf\{\sum_{i=1}^\infty x_{i,r_i}^ah(r_i): E\subset\cup_{i=1}^\infty B(x_i,r_i),\ws x_i\in \mathbb{R}_+^n, \ws r_i \le \varrho \},\end{equation}
where $x_{i,r_i}^a:=\max \{t^a : t\in \overline{B(x_i,r_i)}\}$.\\
Since $\Lambda^{\varrho}_{h,a}(E)$ is a decreasing function of $\varrho$, we can define the (finite or infinite) $a$-Hausdorff measure of $E$ as
\begin{equation}\Lambda_{h,a}(E)=\lim_{\varrho \to 0}\Lambda^{\varrho}_{h,a}(E).\end{equation}\end{defi}

\medskip

\begin{remark} {\rm (1) $h(x,r):=x_{r}^ah(r)$ is an increasing function of $r$, but it depends on $x$ as well, that is the $a$-Hausdorff measure of $E$ depends also on the location of $E$.\\
(2) If $K\subset \mathrm{int}\mathbb{R}^n_+$ is compact, then there are constants $c_i=c_i(a,K)$, $i=1,2$ such that $c_1 \Lambda_{h}(K)\le \Lambda_{h,a}(K)\le c_2\Lambda_{h}(K)$.\\
(3) As in the standard case $\Lambda^{\infty}_{h,a}(E)=0$ if and only if $\Lambda_{h,a}(E)=0$. Of course, for all $\varrho>0$ $\Lambda^{\infty}_{h,a}(E)\le \Lambda^{\varrho}_{h,a}(E)$, and so $\Lambda^{\infty}_{h,a}(E)\le \Lambda_{h,a}(E)$. Oppositely, by standard arguments if $\Lambda^{\infty}_{h,a}(E)>0$, then there is a constant $c$ such that $\Lambda^{\infty}_{h,a}(E)>c>0$. Let $\varrho$ be so small such that $\Lambda^{\varrho}_{h,a}(E)>c$. The for all $\varrho$-covering of $E$, $\cup_{i=1}^\infty B(x_i,r_i)$, $\sum_{i=1}^\infty x_{i,r_i}^ah(r_i)>c$. If $\cup_{j=1}^\infty B(u_j,r_j)$ is not a $\varrho$-covering of $E$, there exists an $r_l>\varrho$, and because $u_j\in \mathbb{R}^n_+$, $\sum_{j=1}^\infty u_{j,r_j}^ah(r_j)>u_{l,r_l}^ah(\varrho)>c(n)\varrho^{|a|}h(\varrho)$. Thus $\Lambda^{\infty}_{h,a}(E)>\min\{c,c(n)\varrho^{|a|}h(\varrho)\}>0$.}
\end{remark}

\medskip

Let us denote by $\mathcal{Q}_{k}$ the set of the dyadic cubes in $\mathbb{R}_+^n$ with edge length $\frac{1}{2^k}$, $k\in\mathbb{Z}$.

\begin{theorem}\label{tH} With the notation above, let $h$ be an increasing function, $E\subset \mathbb{R}_+^n$ and $\mu \in \mathcal{M}(E)$  such that for all balls $\mu(B(x,r))\le h(r)$. Then
$$\mu_a(E)\le \Lambda^{\infty}_{h,a}(E).$$
Let $h$ be an increasing function, $E \subset Q\in \mathcal{Q}_{k}$. Then there is a constant $c$ depending only on $n$, $k$ and $a$ and a measure $\mu \in \mathcal{M}(E)$   satisfying that $\mu(B(x,r))\le h(r)$ for all balls, such that
$$\Lambda^{\infty}_{h,a}(E)\le c\mu_a(E).$$
\end{theorem}

\proof Obviously if $E\subset \cup_{i=1}^\infty B(x_i,r_i)$, then $$\mu_a(E)\le \sum_{i=1}^\infty \mu_a(B(x_i,r_i))\le \sum_{i=1}^\infty x_{i,r_i}^ah(r_i).$$
The first part of the second statement is proved in \cite[page 137]{ah}, namely there are measures $\mu_{ll}$ such that $\mathrm{supp}\mu_{ll}=\{\cup_jQ_j: Q_j\in\mathcal{Q}_l, \ws Q_j\cap E \neq\emptyset\}$ and $\mu_{ll}(Q_i)\le h(r_i)$ for all $Q_i\in\mathcal{Q}_i$, $i=0, \dots , l$, where $r_i=\frac{1}{2^i}$. Moreover $\mu_{ll}$ has constant density on each $Q_j\in \mathcal{Q}_l$. Finally $\mu$ is defined as a weak accumulation point of $\{\mu_{ll}\}$. Then $\mathrm{supp}\mu =E$ and $\mu(Q_k)\le 3^nh(r_k)$ for all $Q_k\in \mathcal{Q}_{k}$, $k\in\mathbb{N}$. It is also pointed out that $E$ has a disjoint covering with dyadic cubes, $E\subset\cup_j Q_j$, such that $Q_j\in\mathcal{Q}_{n_j}$ and $\mu_{ll}(Q_j)=h(r_{n_j})$, $j=1,2,\dots$. Thus $\mu_{ll}(Q)=\sum_j h(r_{n_j})$. Then, with $m_{n_j}^a=\max_{x\in\overline{Q}_j}x^a$, we have
$$\mu_{a,ll}(Q)=\int_Qx^ad\mu_{ll}(x)\ge c\sum_{j, n_j\le l}m_{n_j}^ah(r_{n_j})\ge c \inf \sum_i m_{n_i}^ah(r_{n_i}),$$
where $c=c(n,a,k)$ and the infimum is taken over all finite or denumerable coverings of $E$. So
$$\mu_a(Q)=\mu_a(E)\ge c\inf \sum_im_{n_i}^ah(r_{n_i}).$$
Taking into consideration that a $Q_i \in\mathcal{Q}_i$ can be covered by $c(n)$ balls of radius $r_i$,
$$\Lambda^{\infty}_{h,a}(E)\le c(n)\inf\sum_i x_{i,r_i}^ah(r_i)\le c\inf \sum_i m_{n_i}^ah(r_{n_i})\le c \mu_a(E),$$
where $c=c(n,a,k)$.

\subsection{Capacity of Cantor sets with Bessel external field}

Let $L:=\{l_k\}_{k=0}^\infty$ be a decreasing sequence such that $0<2l_{k+1}<l_k$ for $k\in \mathbb{N}$. Let $C_0$ be a closed interval of length $l_0$. $C_1$ is obtained by removing an open interval of length $l_0-2l_1$ in the middle of $C_0$, etc., $C_k$ consists of $2^k$ closed intervals of length $l_k$. Let $C_k^n:=C_k\times\dots \times C_k$, the Cartesian product of $n$ copies of $C_k$. Let $C_L:=\cap_{k=0}^\infty C_k^n$. $C_L=C_L(n,Q)$, where $Q=C_0\times \dots \times C_0$, the cube which contains $C_L$.

\medskip

\note

Let $C_k^n=C_k^n(Q,L)=\cup_{i=1}^{2^{nk}}q_{k,i}$ as above, where $q_{k,i}$ are the closed cubes in $C_k^n$ of edge length $l_k$. Let $v_{k,i}^a:=\max_{x \in q_{k,i}} x^a$. Let us denote by
\begin{equation}\label{hlk} h_L(l_k):=h_{Q,L,a}(l_k)=\frac{1}{\sum_{i=1}^{2^{nk}}v_{k,i}^a}.\end{equation}
Obviously, $h(l_k)>h(l_{k+1})$. Let $h_L(r):=h_{Q,L,a}(r)$ be an increasing function on $[0,\infty)$, $h_L(0)=0$ and $h_L(l_k)$ is given by \eqref{hlk}.

\medskip

\begin{theorem}\label{Tcp} Let $0<p\nu<n+|a|$, $C_L(n,Q)$, $h_L=h_{Q,L,a}$ as above. Then $C_{a,\nu,p}(C_L(n,Q))>0$ if and only if
$$\int_0^1\left(\frac{h_L(r)}{r^{n-p\nu}}\right)^{p'-1}\frac{dr}{r}<\infty.$$
\end{theorem}

\proof With the notation above $C_L$ can be covered by $2^{kn}$ balls of radius $l_k \frac{\sqrt{n}}{2}$, $\mathcal{A}(r)\le 2^{kn}$, and if $l_{k+1} \frac{\sqrt{n}}{2}\le r \le l_k \frac{\sqrt{n}}{2}$,
$$\mathcal{B}(r)\le \frac{1}{h_L(l_k)}.$$
Comparing with Theorem \ref{tB} it shows that $C_{a,\nu,p}(C_L)=0$ if $\int_0^1\left(\frac{h_L(r)}{r^{n-p\nu}}\right)^{p'-1}\frac{dr}{r}$ diverges.

On the other hand, considering $h_L$ let us construct the measure $\mu_L$ ensured by Theorem \ref{tH}. In view of Lemma \ref{tk}
$$\chi_{B_+(0,r)}*_a\mu_L(x)\le c r^{|a|}\mu_L(B(x,r))\le cr^{|a|}h_L(r).$$
According to Theorem \ref{cw}
$$\|G_{a,\nu}*_a\mu_L\|_{p',a}^{p'}\le c \int_{\mathbb{R}_+}\int_0^1\left(\frac{\chi_{B_+(0,r)}*_a\mu_L(x)}{r^{n+|a|-p\nu}}\right)^{p'-1}\frac{dr}{r}d\mu_{L,a}(x)$$ $$\le c\int_{\mathbb{R}_+}\int_0^1\left(\frac{r^{|a|}h_L(r)}{r^{n+|a|-p\nu}}\right)^{p'-1}\frac{dr}{r}d\mu_{L,a}(x).$$
In view of \eqref{dc}
$$C_{a,\nu,p}^{\frac{1}{p}}(C_L)\ge \frac{\mu_L(C_L)}{\|G_{a,\nu}*_a\mu_L\|_{p',a}}\ge c \frac{\mu_L(C_L)^{1-\frac{1}{p'}}}{I(h_L)^{\frac{1}{p'}}}$$
which proves the converse statement.

\medskip

\begin{theorem}\label{num} With the notation above and supposing that $l_0=1$ we have that $C_{a,\nu,p}(C_L)>0$ if and only if
$$\sum_{k=0}^\infty \left(l_k^{n-p\nu}\sum_{i=1}^{2^{nk}}v_{k,i}^a\right)^{1-p'}<\infty.$$
\end{theorem}

\proof
First we observe that
\begin{equation}\label{eh}\frac{1}{h_L(l_{k+1})}=\sum_{i=1}^{2^{n(k+1)}}v_{k+1,i}^a=\sum_{j=1}^{2^{nk}}\sum_{i: v_{k+1,i}\in q_j} v_{k+1,i}^a\le 2^n \sum_{i=1}^{2^{nk}}v_{k,i}^a=2^n\frac{1}{h_L(l_{k})}.\end{equation}
In view of \eqref{eh}
$$I(h_L)=\sum_{k=0}^\infty \int_{l_{k+1}}^{l_k}\left(\frac{h_L(r)}{r^{n-p\nu}}\right)^{p'-1}\frac{dr}{r}\le \sum_{k=0}^\infty h_L^{p'-1}(l_k)\int_{l_{k+1}}^{l_k}\left(\frac{1}{r^{n-p\nu}}\right)^{p'-1}\frac{dr}{r}$$ $$\le c\sum_{k=0}^\infty h_L^{p'-1}(l_k)l_{k+1}^{(p\nu-n)(p'-1)}\le c 2^n\sum_{k=0}^\infty h_L^{p'-1}(l_{k+1})l_{k+1}^{(p\nu-n)(p'-1)}$$ $$\le c 2^n\sum_{k=0}^\infty h_L^{p'-1}(l_{k})l_{k}^{(p\nu-n)(p'-1)}.$$
On the other hand
$$I(h_L)\ge c\sum_{k=0}^\infty h_L^{p'-1}(l_{k+1})\left(l_{k+1}^{(p\nu-n)(p'-1)}-l_k^{(p\nu-n)(p'-1)}\right)$$ $$\ge c\sum_{k=0}^\infty h_L^{p'-1}(l_{k+1})l_k^{(p\nu-n)(p'-1)}\ge \frac{c}{2^n}\sum_{k=0}^\infty h_L^{p'-1}(l_{k})l_k^{(p\nu-n)(p'-1)},$$
where in th last but one inequality we used that $2l_{k+1}<l_k$, and then \eqref{eh}.

\medskip

\begin{cor}With the notation above $C_{a,\nu,p}(C_L)>0$ if and only if
$$\sum_{k=0}^\infty \left(l_k^{n-p\nu}2^{nk}\right)^{1-p'}<\infty.$$
\end{cor}

\proof $S_k:=\frac{1}{2^{nk}}\sum_{i=1}^{2^{nk}}v_{k,i}^a$. If $C_L \subset \mathrm{int}\mathbb{R}^n_+$, then $S_k$ obviously can be estimated by a constant from above and from below. If $C_L\subset [0,1]^n$, then
$$1\ge S_k \ge \frac{1}{2^{nk}}\sum_{i: v_{k,i}\in [1-l_1,1]^n}v_{k,i}^a\ge \frac{1}{2^{n}}(1-l_1)^{|a|}.$$
The computation is similar if $C_L\subset Q\not=[0,1]^n$, but $C_L\cap \partial \mathbb{R}^n_+\not= \emptyset$.

\medskip

\begin{proposition}\label{hob}
{\rm (1)} Let $K \subset \mathbb{R}^n_+$ be an arbitrary bounded set and $\varrho>0$.\\
If $\liminf_{r\to 0}\frac{h(r)}{r^n}=0$, then $\Lambda_{h,a}^\varrho(K)=0$.\\
If $\liminf_{r\to 0}\frac{h(r)}{r^n}>0$, then there is a function $\tilde{h}(r)$, increasing and $\tilde{h}(0)=0$ such that $\frac{\tilde{h}(r)}{r^n}$ is decreasing and $\Lambda_{h,a}^\varrho(K)\sim \Lambda_{\tilde{h},a}^\varrho(K)$.\\
{\rm (2)} Let $L:=\{l_k\}$ such that $2l_{k+1}<l_{k}$ and $C_L$ be the corresponding Cantor set. Let $h(r)$ be increasing on $[0,\infty)$, $h(0)=0$.\\
If
$\liminf_{k\to\infty}\frac{h(l_k)}{h_L(l_k)}>0$, then $\Lambda_{h,a}(C_L)>0$.\\
{\rm (3)} With the notation above, there is a constant $c=c(a,Q, n)$ such that
$$\Lambda_{h,a}(C_L)\le c \liminf_{k\to\infty}\frac{h(l_k)}{h_L(l_k)}.$$
\end{proposition}

\proof (1) For any bounded set $\Lambda_{h,a}^\varrho(K)\le c(a,K)\Lambda_{h}^\varrho(K)$, thus \cite[Proposition 5.1.8 (a)]{ah} implies the first statement. To prove the second statement define $\tilde{h}(r)$ with
$$\frac{\tilde{h}(r)}{r^b}:=\inf_{0<t\le r}\frac{h(t)}{t^b}.$$
If $\liminf_{r\to 0}\frac{h(r)}{r^b}>0$, $\tilde{h}(r)>0$ $\forall r>0$. $\frac{\tilde{h}(r)}{r^b}$ decreasing and repeated the chain of ideas of \cite{ah} for an arbitrary $\varepsilon>0$ we choose a $t\in [r,R]$ such that $\frac{h(t)}{t^b}\le (1+\varepsilon)\frac{\tilde{h}(R)}{R^b}\le (1+\varepsilon)\frac{\tilde{h}(t)}{t^b}$ Thus for all $\varepsilon>0$, $\tilde{h}\le h(r)\le h(t)\le  (1+\varepsilon)\tilde{h}(R)\left(\frac{t}{R}\right)^b\le (1+\varepsilon)\tilde{h}(R)$, because $b$ must be positive. That is
$\tilde{h}$ is increasing.\\
Since $\tilde{h}\le h(r)$, it is enough to show that $\Lambda_{h,a}(K)\le c\Lambda_{\tilde{h},a}(K)$. To prove this we assume that $b=n$. Let $K\subset\cup_i B( x_i,r_i)$, $r_i\le \varrho$ is a covering such that $\sum_ix_{i,r_i}^a\tilde{h}(r_i)<\Lambda_{\tilde{h},a}^\varrho(K)+\varepsilon$, where $\varepsilon>0$ is arbitrary. All $B( x_i,r_i)$ can be covered by $c(n)\left(\frac{r_i}{r}\right)^n$ balls of radius $r\le r_i$, $B( x_i,r_i)\subset \cup_jB(x_{i,j},r)$. Taking into account that
$$c(n) \sum_{i=1}^\infty \frac{h(r)}{r^{n}} r_i^n\frac{1}{c(n)}\left(\frac{r}{r_i}\right)^n\sum_{j=1}^{c(n)\left(\frac{r_i}{r}\right)^n}x_{i,j,r}^a \le c(n)\sum_{i=1}^\infty \frac{h(r)}{r^{n}} r_i^n x_{i,r_i}^a,$$ we have
$$\Lambda_{h,a}(K)\le c(n) \sum_{i=1}^\infty \inf_{0< r\le r_i}\frac{h(r)}{r^{n}} r_i^n x_{i,r_i}^a\le c(n)\sum_{i=1}^\infty x_{i,r_i}^a\tilde{h}(r_i),$$
which implies the inequality by choice of the covering.\\
(2)  There is a measure $c(n)\mu_L$ such that $c(n)\mu_L(B(x,r))\le h(r)$ for all $r\le r_0$, cf. \cite[Theorem 5.3.1]{ah}. Thus Theorem \ref{tH} ensures the result.\\
(3) $C_L\subset Q$ can be covered by $c(n)2^{kn}$ balls of radius $l_k$, $B(x_i,l_k)$. As $x_{i,l_k}^a\le cv_{i,k}^a$, if $\varrho\ge l_k$ $(c=c(Q,n,a)$,
$$\Lambda_{\tilde{h},a}^\varrho(C_L)\le c c(n)\sum_{i=1}^{2^{kn}}v_{i,k}^ah(l_k)=c(cn)\frac{h(l_k)}{h_L(l_k)},$$
which implies the statement.

\medskip

Another corollary of Theorem \ref{num} is given below.

\begin{cor}\label{T9} {\rm (a)} Let $h(r)$ be increasing on $[0,\infty)$, $h(0)=0$. Let $0<p\nu\le n$. If
$$\int_0^1\left(\frac{h(r)}{r^{n-p\nu}}\right)^{p'-1}\frac{dr}{r}=\infty,$$
then there exists a compact set $K \subset \mathbb{R}^n_+$ such that $\Lambda_{h,a}(K)>0$ and $C_{a,\nu,p}(K)=0$.\\
{\rm (b)} With the notation above, if
\begin{equation}\label{hpc}\liminf_{r\to 0}\frac{h(r)}{r^{n-p\nu}}=0,\end{equation}
then there exists a compact set $K \subset \mathbb{R}^n_+$ such that $\Lambda_{h,a}(K)=0$ and $C_{a,\nu,p}(K)>0$.
\end{cor}

\proof Comparing Theorem \ref{num} with \cite[Theorem 5.3.2]{ah}, it can be seen that if $0<p\nu\le n$, $C_{a,\nu,p}(C_L)>0$ if and only if $C_{\nu,p}(C_L)>0$.
If $C_L\subset [0,1]^n$, $\Lambda_{h}(C_L)\ge \Lambda_{h,a}(C_L) \ge \Lambda_{h,a}(C_L^1)\ge (1-l_1)^{|a|}\Lambda_{h}(C_L^1)\ge c(l_1, a,n)\Lambda_{h}(C_L)$, where $C_L^1$ is the Cantor set associated with $L^1:=\{l_i\}_{i=1}^\infty$ and located in $[1-l_1,1]^n$. According to \cite[Theorem 5.4.2]{ah} assumptions of (a) ensures that there is a Cantor set $C_L$ with $\Lambda_{h}(C_L)>0$ and $C_{\nu,p}(C_L)=0$, which example fulfils the requirements of part (a). According to \cite[Theorem 5.4.1]{ah} assumptions of part (b) imply that there is a Cantor set $C_L$ with $\Lambda_{h}(C_L)=0$ and $C_{\nu,p}(C_L)>0$, which example fulfils the requirements of part (b).

\medskip

{\bf Construction.}
 Comparing Theorem \ref{Tcp}, Proposition \ref{hob} and Corollary \ref{T9} there is a sequence $L=\{l_k\}_{k=0}^\infty$ with a corresponding Cantor set such that $h(l_k)\sim h_L(l_k)$. If $\frac{h(r)r^{|a|}}{r^n}$ is decreasing, it is not difficult to construct a sequence which generates a Cantor set and fulfils that $h(l_k)=h_L(l_k)$.\\
Indeed, Let $Q=[0,1]^n$. We can assume that $h(1)=1=l_0$. The right endpoints of the intervals of $C_1$ are $l$, $1$. Thus all the coordinates of $v_1,i(l)$ are $l$ or $1$ and so $u_1(l):=\sum_{i=1}^{2^n}v_1,i^a(l)$ is increasing in $l$, and $f_l(l):=\frac{1}{u_1(l)}$ is decreasing, positive $f_1(0)=1$ and $f_1(1)=\frac{1}{2^n}$. So $l_1$ is defined by $h(l_1)=f_1(l_1)$.\\
$u_2(l)= u_1(l_1)+s_2(l)=\sum_{i=1}^{2^{2n}}v_{2,i}^a(l_1,l)$, where $s_2(l)$ is an increasing function of $l$ since the coordinates of $v_{2,i}$ contains the right endpoints of $C_2=C_2(l_1,l)$. $u_2(l_1)=2^nu_1(l_1)$, because the right endpoints $l; l_1; 1-l_1+l; 1$ become $l_1;l_1; 1; 1$. Thus $f_2(l):=\frac{1}{u_2(l)}$ is decreasing, $f_2(l_1)=\frac{1}{2^n}f_1(l_1)$, $f_2(0)>0$, so $l_2$ is defined by $h(l_2)=f_2(l_2)$. ... $u_k(l)=u_{k-1}(l_{k-1})+s_k(l)$, $f_k(l):=\frac{1}{u_k(l)}$ is decreasing, $f_k(l_{k-1})=\frac{1}{2^n}f_{k-1}(l_{k-1})=\frac{1}{2^n}h(l_{k-1})$, $f_k(0)>0$, which defines $l_k$. $L:=\{l_k\}$ is obviously decreasing and $h(l_{k-1})\le 2^nh(l_{k})$. It remains to prove that $L$ defines a Cantor set. We have to show that
$$f_k\left(\frac{l_{k-1}}{2}\right)\le h\left(\frac{l_{k-1}}{2}\right).$$
Consider $u_k\left(\frac{l}{2}\right)=u_{k-1}(l_{k-1})+s_k\left(\frac{l}{2}\right)$, and the members of $s_k\left(\frac{l}{2}\right)$ consist of products of terms $\left( d+\frac{l}{2}\right)^{a_i}$, where $d\ge 0$. Thus $u_k\left(\frac{l}{2}\right)\ge \frac{1}{2^{|a|}}u_k(l)$. So
$$\frac{1}{2^{|a|}}f_k\left(\frac{l_{k-1}}{2}\right)\le f_k(l_{k-1})=\frac{1}{2^n}f_{k-1}(l_{k-1})=\frac{1}{2^n}h(l_{k-1})\le \frac{1}{2^{|a|}}h\left(\frac{l_{k-1}}{2}\right),$$
where in the last inequality we used the assumption.

Notice, that together with Theorem \ref{tB} this leads to the construction of a Cantor-type set in $[0,1]^n$ with "prescribed" B-$p$-capacity. More precisely, if $0< |a|<p\nu$, let $h(r)=r^{n-b}$, where $0<p\nu-b<p\nu-|a|$. Then $0<C_{a,\nu,p}(C_L)<c(n,a)p(p\nu-b)$, where $h(l_k)=h_L(l_k)$.

\medskip

\noindent \small{Department of Analysis, Institute of Mathematics,\newline
Budapest University of Technology and Economics \newline
 M\H uegyetem rkp. 3., H-1111 Budapest, Hungary.
 
 \vspace{3mm}
 
\noindent Department of Analysis, \newline
Alfr\'ed R\'enyi Institute of Mathematics, \newline
Re\'altanoda street 13-15, H-1053, Budapest, Hungary

\medskip

 g.horvath.agota@renyi.hu}

\medskip


\begin{thebibliography}{99}

\bibitem {a} D. R. Adams, Weighted nonlinear potential theory, {\it Trans. of the Amer. Math. Soc.} {\bf 297} (1) (1986) 73-94.

\bibitem {ah} D. R. Adams, L. I. Hedberg, {\it Function Spaces and Potential Theory}, Springer-Verlag Berlin Heidelberg (1996)

\bibitem {aclm} N. Arcozzi, N. Chalmoukis, M. Levi, P. Mozolyako, Two-weight dyadic Hardy's inequalities, (2022) arXiv:2110.05450

\bibitem {abdh} V. Arestov, A. Babenko, M. Deikalova, \'A. Horv\'ath, Nikol'skii inequality between the uniform norm and the integral norm with Bessel weight for entire functions of exponential type on the half-line, {\it Analysis Math.}, {\bf 44} (2018) 21-42.

\bibitem {d} J. Delsarte, Sur une extension de la formule de Taylor, {\it Journ. de Math. pures et appliquees}  {\bf 17} 3 (1938), 213-231.

\bibitem {de} J. Deny, Les potentiels d'energie finie, {\it Acta Math.} {\bf 82} (1950), 107-183.

\bibitem {dls} A. Dzhabrailov, Yu. Luchko, E. Shishkina, Two forms of an inverse operator to the generalized Bessel potential, {it Axioms} {\bf 10}, (3) 10.3390/axioms10030232 (2021).

\bibitem {g} V. S. Guliev, On maximal function and fractional integral, associated with the Bessel differential operator, {\it Mathematical Inequalities and Applications} {\bf 6} (2003) 317-330.

\bibitem {gu2} V. S. Guliev, J. J. Hasanov, Necessary and sufficient conditions for the boundedness of B-Riesz potential in the B-Morrey spaces, {\it J. Math. Anal. Appl.} {\bf 347} (2008) 113-122.

\bibitem {ghn} P. Gurka, P. Harjulehto, A. Nekvinda, Bessel potential spaces with variable exponent, {\it Math. Ineq. and Appl.} {\bf 10} (3) (2007) 661-676.

\bibitem {h} \'A. P. Horv\'ath, Compactness Criteria via Laguerre and Hankel transformations, {\it J. Math. Anal. Appl.} {\bf 507} (2) (2022)  125852

\bibitem {kk} I. A. Kipriyanov, V. I. Kononenko, Fundamental solutions of B-elliptic equations, {\it Differ.Uravn.} {\bf 3} (1967) 114-129.

\bibitem {ly} L. N. Lyakhov, Fundamental solutions of singular differential equations with a Bessel $D_B$ operator, {\it Proc. Steklov Inst.Math.} {\bf 278} (2012) 139-151.

\bibitem{le} B. M. Levitan, Expansion in Fourier series and integrals with Bessel functions, {\it Uspekhi Mat. Nauk} (in Russian) {\bf 6} ( 1951) 102--143.

\bibitem{mp} G. Mingione, G. Palatucci, Developments and perspectives in nonlinear potential theory, {\it Nonlinear Analysis} {\bf 194} (2020) 111452.

\bibitem{p} S. S. Platonov, Bessel harmonic analysis and approximation of functions on the half-line, {\it Izv. RAN, Ser. Mat.},  {\bf 71} (2007), 149-196 (in Russian); translated in {\it Izv. Math.}, {\bf 71}:5 (2007), 1001-1048.

\bibitem{rs} H. Rafeiro, S. Samko, Characterization of the variable exponent Bessel potential spaces via the Poisson semigroup, {\it J. Math. Anal. Appl.} {\bf 365} (2010) 483-497.

\bibitem {esk} E. Shishkina, I. Ekincioglu, C. Keskin, Generalized Bessel potential and its application to non-homogeneous singular screened Poisson equation, {\it Integral Transforms and Special Functions} {\bf 32} (12) (2021) https://doi.org/10.1080/10652469.2020.1867983

\bibitem {s} E. Shishkina, Mean-value theorem for B-harmonic functions, {\it Lobachevskii Journal of Mathematics} {\bf 43} (6) (2022) 1401-1407.

\bibitem {ss} E. Shishkina, S. Sitnik, {\it Transmutations, Singular and Fractional Differential Equations With Applications to Mathematical Physics}, Academic Press, London (2020).

\bibitem {sat} E. B. Saff, V. Totik, {\it Logarithmic Potentials with External Fields} Springer-Verlag Berlin Heidelberg (1997)

\bibitem {st} E. M. Stein, {\it Singular Integrals and Differentiability of Functions}, Princeton University
Press, Princeton, New Jersey, (1970).

\bibitem {t} M. R. Tupputi, Weighted Inequalities in some potential spaces on the upper half space of $R^{n+1}$, {\it Potential Anal} {\bf 42} (2015) 293-309.

\bibitem {v} I. E. Verbitsky, Nonlinear potentials and trace inequalities, {\it Oper. Theory: Adv. Appl.} {\bf 110} (1999) 323-343.

\bibitem {v1} I. E. Verbitsky, Wolff's inequality for intrinsic nonlinear potentials and quasilinear elliptic equations, {\it Nonlinear Analysis} {\bf 194} (2020) 111516.

\end{thebibliography}
\end{document}